# LÉVY PROCESSES: CAPACITY AND HAUSDORFF DIMENSION

By Davar Khoshnevisan and Yimin Xiao

*University of Utah and Michigan State University*

We use the recently-developed multiparameter theory of additive Lévy processes to establish novel connections between an arbitrary Lévy process $X$ in $\mathbf{R}^d$, and a new class of energy forms and their corresponding capacities. We then apply these connections to solve two long-standing problems in the folklore of the theory of Lévy processes.

First, we compute the Hausdorff dimension of the image $X(G)$ of a nonrandom linear Borel set $G \subset \mathbf{R}_+$, where $X$ is an arbitrary Lévy process in $\mathbf{R}^d$. Our work completes the various earlier efforts of Taylor [*Proc. Cambridge Phil. Soc.* **49** (1953) 31–39], McKean [*Duke Math. J.* **22** (1955) 229–234], Blumenthal and Getoor [*Illinois J. Math.* **4** (1960) 370–375, *J. Math. Mech.* **10** (1961) 493–516], Millar [*Z. Wahrsch. verw. Gebiete* **17** (1971) 53–73], Pruitt [*J. Math. Mech.* **19** (1969) 371–378], Pruitt and Taylor [*Z. Wahrsch. Verw. Gebiete* **12** (1969) 267–289], Hawkes [*Z. Wahrsch. verw. Gebiete* **19** (1971) 90–102, *J. London Math. Soc. (2)* **17** (1978) 567–576, *Probab. Theory Related Fields* **112** (1998) 1–11], Hendricks [*Ann. Math. Stat.* **43** (1972) 690–694, *Ann. Probab.* **1** (1973) 849–853], Kahane [*Publ. Math. Orsay* **(83-02)** (1983) 74–105, *Recent Progress in Fourier Analysis* (1985b) 65–121], Becker-Kern, Meerschaert and Scheffler [*Monatsh. Math.* **14** (2003) 91–101] and Khoshnevisan, Xiao and Zhong [*Ann. Probab.* **31** (2003a) 1097–1141], where $\dim X(G)$ is computed under various conditions on $G$, $X$ or both.

We next solve the following problem [Kahane (1983) *Publ. Math. Orsay* **(83-02)** 74–105]: *When $X$ is an isotropic stable process, what is a necessary and sufficient analytic condition on any two disjoint Borel sets $F, G \subset \mathbf{R}_+$ such that with positive probability, $X(F) \cap X(G)$ is nonempty*? Prior to this article, this was understood only in the case that $X$ is a Brownian motion [Khoshnevisan (1999) *Trans. Amer. Math. Soc.* **351** 2607–2622]. Here, we present a solution to Kahane's problem for an arbitrary Lévy process $X$, provided the distribution of

Received November 2003; revised June 2004.
Supported in part by NSF Grant DMS-01-03939.

$X(t)$ is mutually absolutely continuous with respect to the Lebesgue measure on $\mathbf{R}^d$ for all $t > 0$.

As a third application of these methods, we compute the Hausdorff dimension and capacity of the preimage $X^{-1}(F)$ of a nonrandom Borel set $F \subset \mathbf{R}^d$ under very mild conditions on the process $X$. This completes the work of Hawkes [*Probab. Theory Related Fields* **112** (1998) 1–11] that covers the special case where $X$ is a subordinator.

**1. Introduction.** It has been long known that a typical Lévy process $X := \{X(t)\}_{t \geq 0}$ in $\mathbf{R}^d$ maps a Borel set $G \subset \mathbf{R}_+$ to a random fractal $X(G)$. For example, Blumenthal and Getoor (1960) have demonstrated that when $X$ is an $\alpha$-stable Lévy process in $\mathbf{R}^d$, then for all Borel sets $G \subset \mathbf{R}_+$,

$$(1.1) \qquad \dim X(G) = d \wedge \alpha \dim G \qquad \text{a.s.,}$$

where dim denotes Hausdorff dimension. In plain words, an $\alpha$-stable process maps a set of Hausdorff dimension $\beta$ to a set of Hausdorff dimension $d \wedge \alpha\beta$. For earlier works in this area, see Taylor (1953) and McKean (1955), and for background on Hausdorff dimension and its properties, see Falconer (1990) and Mattila (1995).

Blumenthal and Getoor (1961) extended (1.1) to a broad class of Lévy processes. For this purpose, they introduced the upper index $\beta$ and lower indices $\beta', \beta''$ of a general Lévy process $X$ and, in addition, the lower index $\sigma$ of a subordinator. Blumenthal and Getoor [(1961), Theorems 8.1 and 8.5] established the following upper and lower bounds for $\dim X(G)$ in terms of the upper index $\beta$ and lower indices $\beta'$ and $\beta''$ of $X$: For every $G \subset \mathbf{R}_+$, almost surely,

$$(1.2) \qquad \begin{aligned} &\dim X(G) \leq \beta \dim G \qquad \text{if } \beta < 1, \\ &\dim X(G) \geq \begin{cases} \beta' \dim G, & \text{if } \beta' \leq d, \\ 1 \wedge \beta'' \dim G, & \text{if } \beta' > d = 1. \end{cases} \end{aligned}$$

They showed, in addition, that when $X$ is a subordinator, then

$$(1.3) \qquad \sigma \dim G \leq \dim X(G) \leq \beta \dim G \qquad \text{a.s.}$$

The restriction $\beta < 1$ of (1.2) was removed subsequently by Millar [(1971), Theorem 5.1]. Blumenthal and Getoor [(1961), page 512] conjectured that, given a Borel set $G \subset [0,1]$, there exists a constant $\lambda(X, G)$ such that

$$(1.4) \qquad \dim X(G) = \lambda(X, G) \qquad \text{a.s.}$$

Moreover, they asked a question that we rephrase as follows: *Given a Lévy process $X$, is it always the case that $\dim X(G) = \dim X([0,1]) \cdot \dim G$ for all nonrandom Borel sets $G \subseteq \mathbf{R}_+$?* Surprisingly, the answer to this question is "no" [Hendricks (1972) and Hawkes and Pruitt (1974)]. To paraphrase from Hawkes and Pruitt [(1974), page 285], in general, $\dim X(G)$ depends on other characteristics of the set $G$ than its Hausdorff dimension. Except in the case where $X$ is a subordinator [Hawkes (1978), Theorem 3], this question had remained unanswered.

One of our original aims was to identify precisely what these characteristics are. As it turns out, the complete answer is quite unusual; see Theorem 2.2. For an instructive example, also Theorem 7.1.

In the slightly more restrictive case that $X$ is a *symmetric* $\alpha$-stable Lévy process, Kahane [(1985b), see Theorem 8] proved that for any Borel set $G \subset \mathbf{R}_+$,

$$(1.5) \qquad \mathcal{H}_\gamma(G) = 0 \implies \mathcal{H}_{\alpha\gamma}(X(G)) = 0 \qquad \text{a.s.}$$

Here, $\mathcal{H}_\beta$ denotes the $\beta$-dimensional Hausdorff measure. If, in addition, we assume that $\alpha\gamma < d$, then Kahane's theorem states further that

$$(1.6) \qquad \mathcal{C}_\gamma(G) > 0 \implies \mathcal{C}_{\alpha\gamma}(X(G)) > 0 \qquad \text{a.s.,}$$

where $\mathcal{C}_\beta$ denotes the $\beta$-dimensional Bessel–Riesz capacity which we recall at the end of this introduction.

As regards a converse to (1.6), Hawkes (1998) has recently proven that if $X$ is a stable *subordinator* of index $\alpha \in (0,1)$, then for any Borel set $G \subset \mathbf{R}_+$ and for all $\gamma \in (0,1)$,

$$(1.7) \qquad \mathcal{C}_\gamma(G) > 0 \iff \mathcal{C}_{\alpha\gamma}(X(G)) > 0 \qquad \text{a.s.}$$

The arguments devised by Hawkes (1998) use specific properties of stable subordinators, and do not apply to other stable processes. On the other hand, Kahane's proof of (1.5) depends crucially on the self-similarity of strictly stable processes. Thus, these methods do not apply to more general Lévy processes.

Our initial interest in such problems came from the surprising fact that the existing literature does not seem to have a definitive answer for the following question:

QUESTION 1.1. Can one find a nontrivial characterization of when $\mathcal{C}_\gamma(X(G))$ is positive for a $d$-dimensional Brownian motion $X$?

The main purpose of this paper is to close the gaps in (1.5) and (1.6) and their counterparts for the preimages of $X$. While doing so, we also answer Question 1.1 in the affirmative. [The answer is the most natural one: "$\mathcal{C}_\gamma(X(G)) > 0$ if and only if $\mathcal{C}_{\gamma/2}(G) > 0$"; cf. Theorem 7.1.]

Our methods rely on a great deal of the recently-developed potential theory for additive Lévy processes; see Khoshnevisan and Xiao (2002, 2003) and Khoshnevisan, Xiao and Zhong (2003a). While the present methods are quite technical, they have the advantage of being adaptable to very general settings. Therefore, instead of working with special processes such as stable processes, we state our results for broad classes of Lévy processes. Moreover, the present methods allow us to solve the following long-standing problem: "Given a Lévy process $X$ in $\mathbf{R}^d$, and two disjoint sets $F, G \subset \mathbf{R}_+$, when is $X(F) \cap X(G)$ nonempty?" Kahane (1983) studied this problem for a symmetric stable Lévy process $X$ in $\mathbf{R}^d$ and proved that

$$
\begin{aligned}
(1.8) \quad \mathcal{C}_{d/\alpha}(F \times G) > 0 &\implies \mathrm{P}\{X(F) \cap X(G) \neq \varnothing\} > 0 \\
&\implies \mathcal{H}_{d/\alpha}(F \times G) > 0.
\end{aligned}
$$

Kahane [(1983), page 90] conjectured that $\mathcal{C}_{d/\alpha}(F \times G) > 0$ is necessary and sufficient for $\mathrm{P}\{X(F) \cap X(G) \neq \varnothing\} > 0$. Until now, this problem had been solved only when $X$ is a Brownian motion [Khoshnevisan (1999), Theorem 8.2].

For a Lévy process $X$ in $\mathbf{R}^d$, we investigate the Hausdorff dimension and capacity of the preimage $X^{-1}(R)$, where $R \subset \mathbf{R}^d$ is a Borel set. When $X$ is *isotropic* $\alpha$-stable, Hawkes (1971) has proven that if $\alpha \geq d$, then for every Borel set $R \subset \mathbf{R}^d$,

$$
(1.9) \qquad \dim X^{-1}(R) = \frac{\alpha + \dim R - d}{\alpha} \qquad \text{a.s.,}
$$

and if $\alpha < d$, then

$$
(1.10) \qquad \sup\{\theta > 0 : \mathrm{P}\{\dim X^{-1}(R) \geq \theta\} > 0\} = \frac{\alpha + \dim R - d}{\alpha}.
$$

More recently, Hawkes (1998) has studied the capacity of $X^{-1}(R)$ further in the case that $X$ is a symmetric $\alpha$-stable Lévy process in $\mathbf{R}$. We are able to extend his result to a large class of Lévy processes; see Theorem 3.1 and Corollary 3.2 below.

We conclude this introduction by introducing some notation that will be used throughout.

We write $\mathcal{P}(F)$ for the collection of all Borel-regular probability measures on a given Borel space $F$.

Given a Borel measurable function $f : \mathbf{R}^d \to [0, \infty]$, we define the "$f$-energy" [of some $\mu \in \mathcal{P}(\mathbf{R}^d)$] and "$f$-capacity" (of some measurable $G \subset \mathbf{R}^d$) as follows:

$$
\begin{aligned}
(1.11) \quad \mathcal{E}_f(\mu) &:= \iint f(x - y) \mu(dx) \mu(dy), \\
\mathcal{C}_f(G) &:= \left[ \inf_{\mu \in \mathcal{P}(G)} \mathcal{E}_f(\mu) \right]^{-1}.
\end{aligned}
$$

We refer to such a function $f$ as a *gauge* function. Occasionally, we write $\mathcal{E}_f(\mu)$ for a bounded measurable $f : \mathbf{R}^d \to \mathbb{C}$, as well.

Given a number $\beta > 0$, we reserve $\mathcal{C}_\beta$ and $\mathcal{E}_\beta$ for $\mathcal{C}_f$ and $\mathcal{E}_f$, respectively, where the gauge function $f$ is $f(t) := \|t\|^{-\beta}$. $\mathcal{C}_\beta$ and $\mathcal{E}_\beta$ are, respectively, the ($\beta$-dimensional) Bessel–Riesz *capacity* and *energy* to which some references were made earlier. More information about the Bessel–Riesz capacity and its connection to fractals can be found in Mattila (1995), Kahane (1985a) and Khoshnevisan (2002). For a lively discussion of the various connections between random fractals, capacity and fractal dimensions, see Taylor (1986).

An important aspect of our proofs involves artificially expanding the parameter space from $\mathbf{R}_+$ to $\mathbf{R}_+^{1+p}$ for an arbitrary positive integer $p$. For this, we introduce some notation that will be used throughout: Any $\mathbf{t} \in \mathbf{R}^{1+p}$ is written as $\mathbf{t} := (t_0, \vec{t})$, where $\vec{t} := (t_1, \ldots, t_p) \in \mathbf{R}^p$. This allows us to extend any $\mu \in \mathcal{P}(\mathbf{R}_+)$ to a probability measure $\boldsymbol{\mu}$ on $\mathbf{R}_+^{1+p}$ as follows:

$$(1.12) \qquad \boldsymbol{\mu}(d\mathbf{t}) := \mu(dt_0) e^{-\sum_{j=1}^p t_j} \, d\vec{t}.$$

Finally, the Lebesgue measure on $\mathbf{R}^d$ is denoted by $\lambda_d$, and for any integer $k$ and all $x, y \in \mathbf{R}^k$, we write $x \prec y$ in place of the statement that $x$ is less than or equal to $y$, coordinatewise; that is, $x_i \leq y_i$ for all $i \leq k$.

**2. The image of a Borel set.** The first main result of this paper is the most general theorem on the Bessel–Riesz capacity of the image $X(G)$ of a Lévy process $X$ in $\mathbf{R}^d$. This, in turn, provides us with a method for computing $\dim X(G)$ for a nonrandom Borel set $G \subset \mathbf{R}_+$. Our computation involves terms that are solely in terms of the Lévy exponent $\Psi$ of $X$ and the set $G$. Hence, Corollary 2.6 verifies the conjecture (1.4) of Blumenthal and Getoor (1961).

Before stating our formula for $\dim X(G)$, we introduce some notation.

Given any $\xi \in \mathbf{R}^d$, we define the function $\chi_\xi : \mathbf{R} \to \mathbb{C}$ as follows:

$$(2.1) \qquad \chi_\xi^\Psi(x) := \chi_\xi(x) := e^{-|x|\Psi(\mathrm{sgn}(x)\xi)} \qquad \forall x \in \mathbf{R}.$$

We will write the more tedious $\chi_\xi^\Psi$ in favor of $\chi_\xi$ only when there are more than one Lévy exponent in the problem at hand and there may be ambiguity as to which Lévy exponent is in question.

Below are some of the elementary properties of this function $\chi_\xi$.

LEMMA 2.1. *For any* $\xi \in \mathbf{R}^d$, $\sup_{x \in \mathbf{R}} |\chi_\xi(x)| \leq 1$. *Moreover, given any* $\mu \in \mathcal{P}(\mathbf{R}_+)$, $\mathcal{E}_{\chi_\xi}(\mu) \geq 0$ *for all* $\xi \in \mathbf{R}^d$. *In particular,* $\mathcal{E}_{\chi_\xi}(\mu) \in [0, 1]$ *is real-valued.*

PROOF. We note that for any $s, t \geq 0$, and for all $\xi \in \mathbf{R}^d$,

$$(2.2) \qquad \mathrm{E}[e^{i\xi \cdot (X(t) - X(s))}] = \chi_\xi(t - s).$$

This shows that $\chi_\xi$ is pointwise bounded in modulus by one. Moreover, by the Fubini–Tonelli theorem, given any $\mu \in \mathcal{P}(\mathbf{R}_+)$, we can integrate the preceding display $[\mu(dt)\mu(ds)]$ to deduce that

$$\mathcal{E}_{\chi_\xi}(\mu) = \mathrm{E}\bigg[\bigg|\int e^{i\xi \cdot X(t)} \mu(dt)\bigg|^2\bigg], \tag{2.3}$$

which completes our proof. $\square$

We are finally ready to present the first main contribution of this paper. The following theorem closes the gaps in (1.6) and (1.7) for a general Lévy process.

THEOREM 2.2. *Suppose $X := \{X(t)\}_{t\geq 0}$ is a Lévy process in $\mathbf{R}^d$, and denote its Lévy exponent by $\Psi$. Then for any Borel set $G \subset \mathbf{R}_+$, and for all $\beta \in (0,d)$,*

$$\begin{aligned}\mathcal{C}_\beta(X(G)) &= 0 \qquad a.s. \\ &\iff \quad \forall \mu \in \mathcal{P}(G) : \int_{\mathbf{R}^d} \mathcal{E}_{\chi_\xi}(\mu) \|\xi\|^{\beta-d}\, d\xi = +\infty.\end{aligned} \tag{2.4}$$

REMARK 2.3. For a closely-related, though different, result, see Khoshnevisan, Xiao and Zhong [(2003a), Theorem 2.1].

The proof of Theorem 2.2 is long, as it requires a good deal of the multiparameter potential theory of additive Lévy process; thus, this proof is deferred to Section 5. In the meantime, the remainder of this section is concerned with describing some of the consequences of Theorem 2.2.

First of all, note that when $X$ is symmetric, $\chi_\xi$ is a positive real function. Thus, we can apply the theorem of Fubini and Tonelli to deduce the following in the symmetric case:

$$\int_{\mathbf{R}^d} \mathcal{E}_{\chi_\xi}(\mu)\|\xi\|^{\beta-d}\, d\xi = \iint \bigg[\int_{\mathbf{R}^d} e^{-|x-y|\Psi(\xi)} \|\xi\|^{\beta-d}\, d\xi\bigg] \mu(dx)\mu(dy). \tag{2.5}$$

In other words, we have the following consequence of Theorem 2.2:

COROLLARY 2.4. *If $X := \{X(t)\}_{t\geq 0}$ is a symmetric Lévy process in $\mathbf{R}^d$ with Lévy exponent $\Psi$, then for any Borel set $G \subset \mathbf{R}_+$, and all $\beta \in (0,d)$,*

$$\mathcal{C}_\beta(X(G)) = 0 \qquad a.s. \iff \mathcal{C}_{f_{d-\beta}}(G) = 0, \tag{2.6}$$

*where*

$$f_\gamma(x) := \int_{\mathbf{R}^d} e^{-|x|\Psi(\xi)} \|\xi\|^{-\gamma}\, d\xi \qquad \forall x \in \mathbf{R}, \gamma \in (0,d). \tag{2.7}$$

REMARK 2.5. We believe that Corollary 2.4 is true quite generally, but have not been successful in proving this. To see the significance of this conjecture, let us assume further that $f_\gamma$ has the property that as $|x|$ tends to zero, $f_\gamma(x) = O(f_\gamma(2x))$. Then, thanks to Corollary 2.4 and a general Frostman theorem [Taylor (1961), Theorem 1], we deduce that for any Borel set $G \subset \mathbf{R}_+$ with finite $f_{d-\beta}^{-1}$-Hausdorff measure, $\mathcal{C}_\beta(X(G)) = 0$ almost surely. In general, we do not know of such conditions in the nonsymmetric case.

Now let us consider the Hausdorff dimension of the image $X(G)$ of any Borel set $G$ under $X$. By the theorem of Frostman [Khoshnevisan (2002), Theorem 2.2.1, Appendix C, and Mattila (1995), Theorem 8.9], given any Borel set $F \subset \mathbf{R}^d$,

$$(2.8) \qquad \dim F := \sup\{\beta \in (0,d) : \mathcal{C}_\beta(F) > 0\}.$$

Thus, Theorem 2.2 allows us to also compute $\dim X(G)$. Namely, we have the following:

COROLLARY 2.6. *Suppose $X := \{X(t)\}_{t \geq 0}$ is a Lévy process in $\mathbf{R}^d$, and denote its Lévy exponent by $\Psi$. Then for any Borel set $G \subset \mathbf{R}_+$,*

$$(2.9) \quad \dim X(G) = \sup\left\{\beta \in (0,d) : \inf_{\mu \in \mathcal{P}(G)} \int_{\mathbf{R}^d} \mathcal{E}_{\chi_\xi}(\mu) \|\xi\|^{\beta-d} d\xi < +\infty\right\} \quad a.s.$$

*In the symmetric case, this is equivalent to the following:*

$$(2.10) \qquad \dim X(G) = \sup\{\beta \in (0,d) : \mathcal{C}_{f_{d-\beta}}(G) > 0\} \quad a.s.,$$

*where $f_\gamma$ is defined in* (2.7).

Corollary 2.6 computes $\dim X(G)$ in terms of the Lévy exponent $\Psi$ of the process $X$. In particular, it verifies the conjecture of Blumenthal and Getoor [(1961), page 512]. However, our formulas are not so easy to use for a given $\Psi$, because they involve an infinite number of computations [one for each measure $\mu \in \mathcal{P}(G)$]. Next, we mention some simple-to-use bounds that are easily derived from Corollary 2.6.

COROLLARY 2.7. *If $X := \{X(t)\}_{t \geq 0}$ is a symmetric Lévy process in $\mathbf{R}^d$, then for any Borel set $G \subset \mathbf{R}_+$, we almost surely have $I(G) \leq \dim X(G) \leq J(G)$, where*

$$(2.11) \quad \begin{aligned} I(G) &:= \sup\left\{\beta \in (0,d) : \limsup_{r \downarrow 0} \frac{\log f_{d-\beta}(r)}{\log(1/r)} < \dim G\right\} \quad and \\ J(G) &:= \inf\left\{\beta \in (0,d) : \liminf_{r \downarrow 0} \frac{\log f_{d-\beta}(r)}{\log(1/r)} > \dim G\right\}. \end{aligned}$$

*In the above, $\inf \varnothing := d$ and $\sup \varnothing := 0$, and $f_\gamma$ is as in* (2.7).

We also mention the following zero–one law. Among other things, it tells us that the a.s.-condition of Theorem 2.2 is sharp.

PROPOSITION 2.8 (Zero–one law). *For any $\beta \in (0,d)$, and for all Borel sets $G \subset \mathbf{R}_+$,*

$$\text{(2.12)} \qquad \mathrm{P}\{\mathcal{C}_\beta(X(G)) > 0\} = 0 \quad or \quad 1.$$

This proposition is a handy consequence of our proof of Theorem 2.2, and its proof is explicitly spelled out in Remark 5.5 below. In the case that $X$ is a subordinator, the reader can find this in Hawkes [(1998), page 9]. We note that Hawkes' proof works for any pure-jump Lévy process.

**3. The preimage of a Borel set.** Let $X := \{X(t)\}_{t \geq 0}$ be a strictly $\alpha$-stable Lévy process in $\mathbf{R}^d$, and let $p_t(x)$ be the density function of $X(t)$. Taylor (1967) proved that

$$\text{(3.1)} \qquad \Gamma := \{x \in \mathbf{R}^d : p_t(x) > 0 \text{ for some } t > 0\}$$

is an open convex cone in $\mathbf{R}^d$ with the origin as its vertex. To further study the structure of $\Gamma$, Taylor (1967) classified strictly stable Lévy processes into two types: $X$ is of *type A* if $p_1(0) > 0$; otherwise it is of *type B*. He proved that when $X$ is of type A, then $p_t(x) > 0$ for all $t > 0$ and $x \in \mathbf{R}^d$. On the other hand, in the case that $X$ is a type B process, Taylor (1967) conjectured that $\Gamma = \{x \in \mathbf{R}^d : p_t(x) > 0 \text{ for all } t > 0\}$; this was later proved by Port and Vitale (1988). By combining the said results, we can conclude that all strictly stable Lévy processes with $\alpha \geq 1$ are of type A.

Now one can extend Hawkes' results (1.9) and (1.10) to all strictly stable Lévy process $X$ of index $\alpha$ in $\mathbf{R}^d$. It follows from Theorem 1 of Kanda (1976) [see also Bertoin (1996), page 61, and Sato (1999), Theorem 42.30] and the arguments of Hawkes (1971) that if $\alpha \geq d$, then for every Borel set $R \subset \mathbf{R}^d$,

$$\text{(3.2)} \qquad \dim X^{-1}(R) = \frac{\alpha + \dim R - d}{\alpha} \qquad \text{a.s.}$$

On the other hand, if $\alpha < d$, then for every Borel set $R \subset \Gamma$,

$$\text{(3.3)} \qquad \|\dim X^{-1}(R)\|_{L^\infty(\mathrm{P})} = \frac{\alpha + \dim R - d}{\alpha},$$

where negative dimension for a set implies that the set is empty.

Hawkes (1998) has made further progress by proving that whenever $X$ is a symmetric $\alpha$-stable process in $\mathbf{R}$, then for all $\beta \in (0,1)$ that satisfy $\alpha + \beta > 1$, and for every Borel $R \subset \mathbf{R}$,

$$\text{(3.4)} \qquad \mathcal{C}_{(\alpha+\beta-1)/\alpha}(X^{-1}(R)) = 0 \qquad \text{a.s.} \quad \Longleftrightarrow \quad \mathcal{C}_\beta(R) = 0.$$

It is an immediate consequence of the Frostman theorem [Khoshnevisan (2002), Theorem 2.2.1, Appendix C] that (3.4) generalizes (3.3). Equation (3.2) also follows from (3.4), Frostman's theorem and recurrence.

In order to go far beyond symmetric stable processes, we can make use of the potential theory of multiparameter Lévy processes. We indicate this connection by proving the following nontrivial generalization of (3.4). For simplicity, we only consider the Lévy processes with $\Gamma = \mathbf{R}^d$.

THEOREM 3.1. *Let $X := \{X(t)\}_{t \geq 0}$ denote a Lévy process in $\mathbf{R}^d$ with Lévy exponent $\Psi$. If $X$ has transition densities $\{p_t\}_{t>0}$ such that for almost all $(t,y) \in \mathbf{R}_+ \times \mathbf{R}^d$, $p_t(y)$ is strictly positive, then for every Borel set $R \subset \mathbf{R}^d$, and all $\gamma \in (0,1)$,*

(3.5)
$$\mathcal{C}_\gamma(X^{-1}(R)) = 0 \qquad a.s.$$
$$\iff \quad \int_{\mathbf{R}^d} |\widehat{\mu}(\xi)|^2 \mathrm{Re}\left(\frac{1}{1 + \Psi^{1-\gamma}(\xi)}\right) d\xi = +\infty \qquad \forall \mu \in \mathcal{P}(R).$$

Theorem 3.1 and Frostman's theorem (2.8) together prove the following:

COROLLARY 3.2. *Let $X := \{X(t)\}_{t \geq 0}$ denote a Lévy process in $\mathbf{R}^d$ with Lévy exponent $\Psi$. If $X$ has strictly positive transition densities, then for every Borel set $R \subset \mathbf{R}^d$,*

(3.6)
$$\|\dim X^{-1}(R)\|_{L^\infty(\mathrm{P})}$$
$$= \sup\left\{\gamma \in (0,1) : \inf_{\mu \in \mathcal{P}(R)} \int_{\mathbf{R}^d} |\widehat{\mu}(\xi)|^2 \mathrm{Re}\left(\frac{1}{1 + \Psi^{1-\gamma}(\xi)}\right) d\xi < +\infty\right\}.$$

Before commencing with our proof of Theorem 3.1, we develop a simple technical result.

Suppose $X := \{X(t)\}_{t \geq 0}$ is a Lévy process in $\mathbf{R}^d$, and suppose that it has transition densities with respect to the Lebesgue measure $\lambda_d$. In other words, we are assuming that there exist (measurable) functions $\{p_t\}_{t \geq 0}$ such that, for all measurable $f : \mathbf{R}^d \to \mathbf{R}_+$ and all $t \geq 0$, $\mathrm{E}[f(X(t))] = \int_{\mathbf{R}^d} f(y) p_t(y)\, dy$.

Next, we consider a $(1-\gamma)$-stable subordinator $\sigma := \{\sigma(t)\}_{t \geq 0}$ that is independent of the process $X$. [Of course, $\gamma$ is necessarily in $(0,1)$.] Let $v_t$ denote the density function of $\sigma(t)$. It is well known that, for every $t > 0$, $v_t(s) > 0$ for all $s > 0$.

LEMMA 3.3. *If $X := \{X(t)\}_{t \geq 0}$ is a Lévy process in $\mathbf{R}^d$ with Lévy exponent $\Psi$ and transition densities $\{p_t\}_{t>0}$, then the subordinated process $X \circ \sigma$ is a Lévy process with Lévy exponent $\Psi^{1-\gamma}$ and transition densities $(t,y) \mapsto \int_0^\infty p_s(y) v_t(s)\, ds$. Moreover, if $p_s(y) > 0$ for almost all $(s,y) \in \mathbf{R}_+ \times \mathbf{R}^d$, then for every $t > 0$, the density of $X(\sigma(t))$ is positive almost everywhere.*

PROOF. Much of this is well known [Sato (1999), Theorem 30.1], and we content ourselves by deriving the transition densities of $X \circ \sigma$. For any measurable $f : \mathbf{R}^d \to \mathbf{R}_+$, and for all $t \geq 0$, $\mathrm{E}[f(X(\sigma(t)))] = \int_{\mathbf{R}^d} f(y) \mathrm{E}[p_{\sigma(t)}(y)] \, dy$. This verifies the formula for the transition densities of $X \circ \sigma$. The final statement of the lemma follows from the well-known fact that $v_t(s) > 0$ for all $s > 0$. □

PROOF OF THEOREM 3.1. As in Lemma 3.3, we let $\sigma$ denote a $(1-\gamma)$-stable subordinator that starts at the origin, and is independent of $X$. Then, it is well known [Hawkes (1971), Lemma 2] that, for any Borel set $B \subset \mathbf{R}_+$,

$$(3.7) \qquad \mathrm{P}\{B \cap \sigma(\mathbf{R}_+) \neq \varnothing\} > 0 \quad \Longleftrightarrow \quad \mathcal{C}_\gamma(B) > 0.$$

By conditioning on $X$, we obtain the following:

$$(3.8) \ \mathrm{P}\{X^{-1}(R) \cap \sigma(\mathbf{R}_+) \neq \varnothing\} > 0 \quad \Longleftrightarrow \quad \mathrm{P}\{\mathcal{C}_\gamma(X^{-1}(R)) > 0\} > 0.$$

Moreover, it is clear that

$$(3.9) \ \mathrm{P}\{X^{-1}(R) \cap \sigma(\mathbf{R}_+) \neq \varnothing\} > 0 \quad \Longleftrightarrow \quad \mathrm{P}\{R \cap X \circ \sigma(\mathbf{R}_+) \neq \varnothing\} > 0.$$

Thus, by Lemma 4.1 below,

$$(3.10) \ \mathrm{P}\{X^{-1}(R) \cap \sigma(\mathbf{R}_+) \neq \varnothing\} > 0 \quad \Longleftrightarrow \quad \mathrm{E}[\lambda_d(X \circ \sigma(\mathbf{R}_+) \ominus R)] > 0.$$

Because $X \circ \sigma$ is a Lévy process in $\mathbf{R}^d$ with exponent $\Psi^{1-\gamma}(\xi)$ (Lemma 3.3), the remainder of the proof follows from Khoshnevisan, Xiao and Zhong [(2003a), Theorem 1.5] that we restate below as Theorem 4.6; see also Bertoin [(1996), page 60]. □

**4. Background on additive Lévy processes.** In this section we rephrase, as well as refine, some of the potential theory of additive Lévy processes that was established in Khoshnevisan and Xiao (2002, 2003) and Khoshnevisan, Xiao and Zhong (2003a). Our emphasis is on how these results are used in order to compute the Hausdorff dimension of various random sets of interest.

A $p$-parameter, $\mathbf{R}^d$-valued, *additive Lévy process* $\vec{X} = \{\vec{X}(\vec{t})\}_{\vec{t} \in \mathbf{R}_+^p}$ is a multiparameter stochastic process that is defined by

$$(4.1) \qquad \vec{X}(\vec{t}) := \sum_{j=1}^{p} X_j(t_j) \qquad \forall \vec{t} = (t_1, \ldots, t_p) \in \mathbf{R}_+^p.$$

Here, $X_1, \ldots, X_p$ denote independent Lévy processes in $\mathbf{R}^d$. Following the notation in Khoshnevisan and Xiao (2002, 2003), we may denote the random field $\vec{X}$ by

$$(4.2) \qquad \vec{X} := X_1 \oplus \cdots \oplus X_p.$$

These additive random fields naturally arise in the analysis of multiparameter processes such as Lévy's sheets and in the studies of intersections of Lévy processes [Khoshnevisan and Xiao (2002)]. [At first sight, the term "additive Lévy" may be redundant. Indeed, historically, the term "additive process" refers to a process with independent increments. Thus, in this sense every Lévy process is additive. However, we feel strongly that our usage of the term "additive process" is more mathematically sound, as can be seen by considering the additive group $\mathfrak{G}$ created by direct-summing cadlag functions $f_1, \ldots, f_p : \mathbf{R}_+ \to \mathbf{R}^d$ to obtain a function $f : \mathbf{R}_+^p \to \mathbf{R}^d$ defined by $f(\vec{t}) := (f_1 \oplus \cdots \oplus f_p)(\vec{t}) = f_1(t_1) + \cdots + f_p(t_p)$. Therefore, if $X_1, \ldots, X_{p+1}$ are independent Lévy processes, then $t_1 \mapsto X_1(t_1) \oplus X_2(\bullet) \oplus \cdots \oplus X_{1+p}(\bullet)$ is a Lévy process on the infinite-dimensional additive group $\mathfrak{G}$.]

For each $\vec{t} \in \mathbf{R}_+^p$, the characteristic function of $\vec{X}(\vec{t})$ is given by

$$(4.3) \qquad \mathrm{E}[e^{i\xi \cdot \vec{X}(\vec{t})}] = e^{-\sum_{j=1}^p t_j \Psi_j(\xi)} := e^{-\vec{t} \cdot \vec{\Psi}(\xi)} \qquad \forall \xi \in \mathbf{R}^d,$$

where $\vec{\Psi}(\xi) := \Psi_1(\xi) \otimes \cdots \otimes \Psi_p(\xi)$, in tensor notation. We will call $\vec{\Psi}(\xi)$ the *characteristic exponent* of the additive Lévy process $\vec{X}$.

Additive Lévy processes have a theory that extends much of the existing theory of Lévy processes. For instance, corresponding to any additive Lévy process $\vec{X}$, there is a *potential measure* $\vec{U}$ that we define as follows: For all measurable sets $F \subset \mathbf{R}^d$,

$$(4.4) \qquad \vec{U}(F) := \mathrm{E}\left[\int_{\mathbf{R}_+^p} e^{-\sum_{j=1}^p s_j} \mathbf{1}_F(\vec{X}(\vec{s}))\, d\vec{s}\right].$$

If $\vec{U}$ is absolutely continuous with respect to the Lebesgue measure $\lambda_d$, its density is called the 1-*potential density* of $\vec{X}$. There is also a notion of transition densities. However, for technical reasons, we sometimes assume more; see Khoshnevisan and Xiao (2002, 2003). Namely, we say that the process $\vec{X}$ is *absolutely continuous* if for each $\vec{t} \in \mathbf{R}_+^p \setminus \partial \mathbf{R}_+^p$, $e^{-\vec{t} \cdot \vec{\Psi}} \in \mathcal{L}^1(\mathbf{R}^d)$. In this case, for all $\vec{t} \in \mathbf{R}_+^p \setminus \partial \mathbf{R}_+^p$, $\vec{X}(\vec{t})$ has a bounded and continuous density function $p(\vec{t}; \bullet)$, which is described by the following formula:

$$(4.5) \qquad p(\vec{t}; x) := (2\pi)^{-d} \int_{\mathbf{R}^d} e^{-i\xi \cdot x - \vec{t} \cdot \vec{\Psi}(\xi)}\, d\xi \qquad \forall x \in \mathbf{R}^d.$$

We remark that when $\vec{X}$ is absolutely continuous, $\vec{U}$ is absolutely continuous and the 1-potential density is $\int_{\mathbf{R}_+^p} p(\vec{s}; \bullet) e^{-\sum_{j=1}^p s_j}\, d\vec{s}$. See Hawkes [(1979), Lemma 2.1] for a necessary and sufficient condition for the existence of a 1-potential density.

When $\vec{X}$ is absolutely continuous, the following function $\Phi$ is well defined, and is called the *gauge function* for $\vec{X}$:

$$(4.6) \qquad \Phi(\vec{s}) := p(|s_1|, \ldots, |s_p|; 0) \qquad \forall \vec{s} \in \mathbf{R}^p.$$

It is clear that $\Phi(\vec{0}) = +\infty$ and, when $X_1, \ldots, X_p$ are symmetric, $\vec{s} \mapsto \Phi(\vec{s})$ is nonincreasing in each $|s_i|$. It is also not too hard to see that $\mathcal{C}_\Phi(\cdot)$ is a natural capacity in the sense of Choquet [Dellacherie and Meyer (1978)].

In order to apply our previous results [Khoshnevisan and Xiao (2002, 2003) and Khoshnevisan, Xiao and Zhong (2003a)], we first extend a one-parameter theorem of Kahane (1972, 1983) to additive Lévy processes in $\mathbf{R}^d$. In the following, we write $A \ominus B := \{x - y : x \in A, y \in B\}$. (Note that when either $A$ or $B$ is the empty set $\varnothing$, then $A \ominus B = \varnothing$.)

LEMMA 4.1. *Let $\vec{X}$ be a p-parameter additive Lévy process in $\mathbf{R}^d$. We assume that, for every $t \in (0,\infty)^p$, the distribution of $\vec{X}(\vec{t})$ is mutually absolutely continuous with respect to $\lambda_d$. Then for all Borel sets $G \subset (0,\infty)^p$ and $F \subset \mathbf{R}^d$, the following are equivalent:*

1. *with positive probability, $G \cap \vec{X}^{-1}(F) \neq \varnothing$;*
2. *with positive probability, $F \cap \vec{X}(G) \neq \varnothing$;*
3. *with positive probability, $\lambda_d(F \ominus \vec{X}(G)) > 0$.*

PROOF. It is clear that $1 \Leftrightarrow 2$. To prove $2 \Leftrightarrow 3$, we note that part 2 is equivalent to the following:

$$(4.7) \qquad \exists \delta > 0 \quad \text{such that} \quad \mathrm{P}\{F \cap \vec{X}(G \cap (\delta,\infty)^p) \neq \varnothing\} > 0.$$

Hence, without loss of generality, we can assume that $G \subset (\delta,\infty)^p$ for some $\delta > 0$. By our assumption, we may choose $\vec{a} \in (0,\infty)^p$ such that:

(i) $\vec{a} \prec \vec{t}$ for all $\vec{t} \in G$.
(ii) The distribution of $\vec{X}(\vec{a})$ is equivalent to $\lambda_d$.

Next, define the additive Lévy process $\vec{X}_{\vec{a}} := \{\vec{X}_{\vec{a}}(\vec{t})\}_{\vec{t} \in \mathbf{R}^p_+}$ by

$$(4.8) \qquad \vec{X}_{\vec{a}}(\vec{t}) := \vec{X}(\vec{t} + \vec{a}) - \vec{X}(\vec{a}) \qquad \forall \vec{t} \in \mathbf{R}^p_+.$$

Then, we point out that

$$(4.9) \qquad F \cap \vec{X}(G) = \varnothing \iff \vec{X}(\vec{a}) \notin F \ominus \vec{X}_{\vec{a}}(G - \vec{a}).$$

Since $\vec{X}(\vec{a})$ is independent of the random Borel set $F \ominus \vec{X}_{\vec{a}}(G - \vec{a})$ and the distribution of $\vec{X}(\vec{a})$ is equivalent to $\lambda_d$, we have

$$(4.10) \quad \begin{aligned} \vec{X}(\vec{a}) &\notin F \ominus \vec{X}_{\vec{a}}(G - \vec{a}) \qquad \text{a.s.} \\ &\iff \lambda_d(F \ominus \vec{X}_{\vec{a}}(G - \vec{a})) = 0 \qquad \text{a.s.} \end{aligned}$$

Note that $\vec{X}_{\vec{a}}(G - \vec{a}) = \vec{X}(G) \ominus \{\vec{X}(\vec{a})\}$, so that the translation invariance of the Lebesgue measure, (4.9) and (4.10) imply that

$$(4.11) \quad F \cap \vec{X}(G) = \varnothing \quad \text{a.s.} \iff \lambda_d(F \ominus \vec{X}(G)) = 0 \quad \text{a.s.}$$

This proves $2 \Leftrightarrow 3$, whence the lemma. □

The following theorem connects the positiveness of the Lebesgue measure of the range $X(G)$ and the hitting probability of the level set $X^{-1}(a)$ to a class of natural capacities. It is a consequence of the results in Khoshnevisan and Xiao [(2002), Theorem 5.1, and (2003)] and Lemma 4.1.

THEOREM 4.2. *Suppose $X_1, \ldots, X_p$ are $p$ independent symmetric Lévy processes on $\mathbf{R}^d$, and let $\vec{X}$ denote $X_1 \oplus \cdots \oplus X_p$, which we assume is absolutely continuous with an a.e.-positive density function at every time $\vec{t} \in \mathbf{R}_+^p$. Let $\Phi$ be the gauge function of $\vec{X}$. Then for every Borel set $G \subset (0, \infty)^p$, the following are equivalent:*

1. $\mathcal{C}_\Phi(G) > 0$.
2. *With positive probability, $\lambda_d(\vec{X}(G)) > 0$.*
3. *Any $a \in \mathbf{R}^d$ can be in the random set $\vec{X}(G)$ with positive probability.*

REMARK 4.3. Theorem 4.2 asserts that, for every Borel set $G \subset (0, \infty)^p$,

$$\mathcal{C}_\Phi(G) > 0 \iff \vec{X}^{-1}(\{0\}) \cap G \neq \varnothing \quad \text{with positive probability.} \tag{4.12}$$

In fact, (4.12) holds even without the assumption that the density function of $\vec{X}(\vec{t})$ is positive almost everywhere; see Corollary 2.13 of Khoshnevisan and Xiao (2002). In Section 6 we apply this minor variation of Theorem 4.2 to derive Theorem 6.5.

We recall that $X_1 \oplus \cdots \oplus X_p$ is *additive $\alpha$-stable* if $X_1, \ldots, X_p$ are independent isotropic $\alpha$-stable processes. The following is a consequence of Lemma 4.1 and Khoshnevisan, Xiao and Zhong [(2003a), Theorem 7.2], which improves the earlier results of Hirsch (1995), Hirsch and Song (1995a, b) and Khoshnevisan (2002).

THEOREM 4.4. *Suppose $\vec{X} := X_1 \oplus \cdots \oplus X_p$ is an additive $\alpha$-stable process in $\mathbf{R}^d$. Then, $\dim \vec{X}(\mathbf{R}_+^p) = \alpha p \wedge d$, a.s. Moreover, for all Borel sets $F \subset \mathbf{R}^d$, the following are equivalent:*

1. $\mathcal{C}_{d-\alpha p}(F) > 0$.
2. *With positive probability, $\lambda_d(F \oplus \vec{X}(\mathbf{R}_+^p)) > 0$.*
3. *$F$ is not polar for $\vec{X}$; that is, with positive probability, $F \cap \vec{X}(\mathbf{R}_+^p \setminus \{\vec{0}\}) \neq \varnothing$.*

REMARK 4.5. Note that the second part of Theorem 4.4 is of interest only in the case that $\alpha p \leq d$. When $\alpha p > d$, $X$ hits every point in $\mathbf{R}^d$ almost surely. Therefore, $\vec{X}(\mathbf{R}_+^p) = \mathbf{R}^d$, a.s. In this case, there is a rich theory of local times and level sets [Khoshnevisan, Xiao and Zhong (2003b)].

PROOF OF THEOREM 4.4. The first statement regarding the dimension of $\vec{X}(\mathbf{R}_+^p)$ follows from Khoshnevisan, Xiao and Zhong [(2003a), Theorem 1.6], whereas $1 \Leftrightarrow 2$ for all compact sets $F$ is precisely Theorem 7.2 of Khoshnevisan, Xiao and Zhong (2003a). In the following we first prove $2 \Leftrightarrow 3$ and then use it to remove the compactness restriction in $1 \Leftrightarrow 2$.

For every $\vec{t} \in \mathbf{R}_+^p \setminus \{\vec{0}\}$, the distribution of $\vec{X}(\vec{t})$ has a strictly positive and continuous density. We write $\mathbf{R}_+^p \setminus \{\vec{0}\} = (0, \infty)^p \cup (\partial \mathbf{R}_+^p \setminus \{\vec{0}\})$. Lemma 4.1 implies that, for every Borel set $F \subset \mathbf{R}^d$,

$$
\begin{aligned}
&\mathrm{P}\{\lambda_d\{F \oplus \vec{X}((0,\infty)^p)\} > 0\} > 0 \\
&\quad \Longleftrightarrow \quad \mathrm{P}\{F \cap \vec{X}((0,\infty)^p) \neq \varnothing\} > 0.
\end{aligned}
\tag{4.13}
$$

For the boundary $\partial \mathbf{R}_+^p \setminus \{\vec{0}\}$, we apply Lemma 4.1 to additive stable processes that have fewer than $p$ parameters to obtain

$$
\mathrm{P}\{\lambda_d\{F \oplus \vec{X}(\partial \mathbf{R}_+^p)\} > 0\} > 0 \quad \Longleftrightarrow \quad \mathrm{P}\{F \cap \vec{X}(\partial \mathbf{R}_+^p \setminus \{\vec{0}\}) \neq \varnothing\} > 0.
\tag{4.14}
$$

Therefore, we have proven $2 \Leftrightarrow 3$ for all Borel sets $F \subset \mathbf{R}^d$. From the above, we derive that, for every compact set $F \subset \mathbf{R}^d$, $1 \Leftrightarrow 3$. But $\mathcal{C}_{d-\alpha p}(\cdot)$ and $\mathsf{C}(\cdot)$ are both Choquet capacities, where

$$
\mathsf{C}(F) = \mathrm{P}\{F \cap \vec{X}(\mathbf{R}_+^p \setminus \{\vec{0}\}) \neq \varnothing\}.
\tag{4.15}
$$

Thus, the compactness restriction on $F$ can be removed by Choquet's capacibility theorem [Dellacherie and Meyer (1978)], whence the validity of $1 \Leftrightarrow 2$ in general. $\square$

We conclude this section by recalling the main results of Khoshnevisan, Xiao and Zhong (2003a). The first is from Khoshnevisan, Xiao and Zhong [(2003a), Theorem 1.5], which can be applied to compute the Hausdorff dimension of the range of an arbitrary Lévy process; for earlier progress on this problem, see Pruitt (1969).

THEOREM 4.6 [Khoshnevisan, Xiao and Zhong (2003a), Theorem 1.5]. *Consider a p-parameter additive Lévy process* $\vec{X} := \{\vec{X}(\vec{t})\}_{\vec{t} \in \mathbf{R}_+^p}$ *in* $\mathbf{R}^d$ *with Lévy exponent* $\Psi$. *Suppose that there exists a constant* $c > 0$ *such that, for all* $\xi \in \mathbf{R}^d$,

$$
\mathrm{Re} \prod_{j=1}^p \{1 + \Psi_j(\xi)\}^{-1} \geq c \prod_{j=1}^p \mathrm{Re}\{1 + \Psi_j(\xi)\}^{-1}.
\tag{4.16}
$$

Then, given a Borel set $F \subset \mathbf{R}^d$, $\mathrm{E}[\lambda_d(\vec{X}(\mathbf{R}_+^p) \oplus F)] > 0$ if and only if there exists $\mu \in \mathcal{P}(F)$ such that

$$(4.17) \qquad \int_{\mathbf{R}^d} |\widehat{\mu}(\xi)|^2 \prod_{j=1}^p \mathrm{Re}\{1 + \Psi_j(\xi)\}^{-1} d\xi < +\infty.$$

As a corollary to this, Khoshnevisan, Xiao and Zhong [(2003a), Theorem 1.6] obtained the following refinement of the results of Pruitt (1969):

COROLLARY 4.7. *If $X$ is a Lévy process in $\mathbf{R}^d$ with Lévy exponent $\Psi$, then a.s.,*

$$(4.18) \quad \begin{aligned} &\dim X(\mathbf{R}_+) \\ &= \sup\left\{\gamma \in (0,d) : \int_{\{\xi \in \mathbf{R}^d : \|\xi\| \geq 1\}} \mathrm{Re}\left(\frac{1}{1 + \Psi(\xi)}\right) \frac{d\xi}{\|\xi\|^{d-\gamma}} < +\infty\right\}. \end{aligned}$$

The next requisite result is from Khoshnevisan, Xiao and Zhong [(2003a), Theorem 2.1 and Lemma 2.4], which characterizes $\mathrm{E}[\lambda_d(\vec{X}(G))] > 0$ completely in terms of its Lévy exponent $\Psi$ and $G$. Notice that it is more general than $1 \Leftrightarrow 2$ in Theorem 4.2.

THEOREM 4.8. *Suppose $\vec{X} := \{\vec{X}(\vec{t})\}_{\vec{t} \in \mathbf{R}_+^p}$ is a $p$-parameter additive Lévy process in $\mathbf{R}^d$ with Lévy exponent $\Psi$. Then, given a Borel set $G \subset \mathbf{R}_+^p$, $\mathrm{E}[\lambda_d(\vec{X}(G))] > 0$ if and only if there exists $\mu \in \mathcal{P}(G)$ such that*

$$\int_{\mathbf{R}^d} \mathcal{E}_{\otimes_{j=1}^p \chi_\xi^{\Psi_j}}(\mu) \, d\xi < +\infty,$$

*where*

$$(4.19) \quad \mathcal{E}_{\otimes_{j=1}^p \chi_\xi^{\Psi_j}}(\mu) = \int_{\mathbf{R}_+^p} \int_{\mathbf{R}_+^p} e^{-\sum_{j=1}^p |s_j - t_j| \Psi_j(\mathrm{sgn}(s_j - t_j)\xi)} \mu(d\vec{s}) \mu(d\vec{t}).$$

**5. Proof of Theorem 2.2.** Next, for any fixed $\alpha \in (0, 2]$, we introduce $p$ independent isotropic $\alpha$-stable Lévy processes $X_1, \ldots, X_p$ in $\mathbf{R}^d$, each of which is normalized as follows:

$$(5.1) \qquad \mathrm{E}[e^{i\xi \cdot X_l(u)}] = e^{-u\|\xi\|^\alpha} \qquad \forall \xi \in \mathbf{R}^d, \ u \geq 0, \ l = 1, \ldots, p.$$

We assume that $X_1, \ldots, X_p$ are independent of the Lévy process $X$ and then consider the *additive Lévy process* $\{A(\mathbf{t})\}_{\mathbf{t} \in \mathbf{R}_+^{1+p}}$; this is the $(1+p)$-parameter random field that is prescribed by the following:

$$(5.2) \qquad A(\mathbf{t}) := X(t_0) + X_1(t_1) + \cdots + X_p(t_p) \qquad \forall \mathbf{t} \in \mathbf{R}_+^{1+p}.$$

For this random field and any $\mu \in \mathcal{P}(\mathbf{R}_+)$, we consider the random measure $\mathrm{O}_\mu$ on $\mathbf{R}^d$ defined by

$$\mathrm{O}_\mu(f) := \int_{\mathbf{R}_+^{1+p}} f(A(\mathbf{t}))\boldsymbol{\mu}(d\mathbf{t}), \tag{5.3}$$

where $\boldsymbol{\mu}$ is defined by (1.12). This is well defined for all nonnegative measurable $f : \mathbf{R}^d \to \mathbf{R}_+$, for instance.

LEMMA 5.1. *For all probability measures $\mu$ on $\mathbf{R}_+$, and $\xi \in \mathbf{R}^d$,*

$$2^{-p}(1+\|\xi\|)^{-\alpha p}\mathcal{E}_{\chi_\xi}(\mu) \leq \mathrm{E}[|\widehat{\mathrm{O}_\mu}(\xi)|^2] \tag{5.4}$$
$$\leq 2^{p\alpha}(1+\|\xi\|)^{-\alpha p}\mathcal{E}_{\chi_\xi}(\mu).$$

PROOF. By Khoshnevisan, Xiao and Zhong [(2003a), Lemma 2.4], for all $\xi \in \mathbf{R}^d$,

$$\mathrm{E}[|\widehat{\mathrm{O}_\mu}(\xi)|^2]$$
$$= \iint_{\mathbf{R}_+^{1+p}\times\mathbf{R}_+^{1+p}} e^{-\sum_{j=1}^p |s_j-t_j|\|\xi\|^\alpha} e^{-|s_0-t_0|\Psi(\operatorname{sgn}(s_0-t_0)\xi)} \boldsymbol{\mu}(d\mathbf{t})\boldsymbol{\mu}(d\mathbf{s}). \tag{5.5}$$

On the other hand, it is easy to see that

$$\iint_{\mathbf{R}_+^p \times \mathbf{R}_+^p} e^{-\sum_{j=1}^p |s_j-t_j|\|\xi\|^\alpha - \sum_{j=1}^p (s_j+t_j)} d\vec{t}\, d\vec{s} = (1+\|\xi\|^\alpha)^{-p}. \tag{5.6}$$

Therefore,

$$\mathrm{E}[|\widehat{\mathrm{O}_\mu}(\xi)|^2] = (1+\|\xi\|^\alpha)^{-p}\mathcal{E}_{\chi_\xi}(\mu). \tag{5.7}$$

To finish the proof, we note merely that

$$\tfrac{1}{2}(1+\|\xi\|^\alpha) \leq (1+\|\xi\|)^\alpha \leq 2^\alpha(1+\|\xi\|^\alpha). \tag{5.8}$$

(For the upper bound, consider the cases $\|\xi\| \leq 1$ and $\|\xi\| > 1$ separately.) □

We obtain the following upon integrating the preceding lemma $[d\xi]$:

LEMMA 5.2. *For all $\mu \in \mathcal{P}(\mathbf{R}_+)$,*

$$2^{-p}(2\pi)^d \mathrm{Q}_\mu^{\alpha p}(\mathbf{R}^d) \leq \mathrm{E}[\|\widehat{\mathrm{O}_\mu}\|_{\mathcal{L}^2(\mathbf{R}^d)}^2] \leq 2^{p\alpha}(2\pi)^d \mathrm{Q}_\mu^{\alpha p}(\mathbf{R}^d), \tag{5.9}$$

*where for any $\beta > 0$,*

$$\frac{\mathrm{Q}_\mu^\beta(d\xi)}{d\xi} := (2\pi)^{-d}\frac{\mathcal{E}_{\chi_\xi}(\mu)}{(1+\|\xi\|)^\beta}. \tag{5.10}$$

REMARK 5.3. Since $\chi_\xi$ is bounded by 1 (Lemma 2.1), for any probability measure $\mu$ on $\mathbf{R}_+$, $\mathcal{E}_{\chi_\xi}(\mu) \leq 1$. Thus, for all $\beta \in (0, d)$,

$$(5.11) \qquad \mathrm{Q}_\mu^\beta(\mathbf{R}^d) < \infty \quad \Longleftrightarrow \quad \int_{\mathbf{R}^d} \mathcal{E}_{\chi_\xi}(\mu) \|\xi\|^{-\beta} \, d\xi < +\infty.$$

Next, we develop a variant of Lemma 5.2. In order to describe it, it is convenient to put all subsequent Lévy processes on the canonical probability space defined by all cadlag paths from $\mathbf{R}_+$ into $\mathbf{R}^d$; see Khoshnevisan, Xiao and Zhong [(2003a), pages 1107 and 1108] for the details of this more-or-less standard construction. Then, we can define the measure $\mathrm{P}_x$, for each $x \in \mathbf{R}^d$, as the measure that starts the process $A$ at $A(0) = x$. Formally speaking, we have $\mathrm{P}_x := \mathrm{P} \circ (A(0) + x)^{-1}$. Since $\lambda_d$ denotes the Lebesgue measure on the Borel subsets of $\mathbf{R}^d$, we can then define

$$(5.12) \qquad \mathrm{P}_{\lambda_d}(W) := \int_{\mathbf{R}^d} \mathrm{P}_x(W) \, dx \quad \text{and} \quad \mathrm{E}_{\lambda_d}[Y] := \int Y \, d\mathrm{P}_{\lambda_d},$$

for all measurable subsets $W$ of the path space and all positive random variables $Y$. An important fact about additive Lévy processes is that they satisfy the Markov property with respect to the $\sigma$-finite measure $\mathrm{P}_{\lambda_d}$. See Khoshnevisan and Xiao [(2002), Proposition 5.8] or Khoshnevisan, Xiao and Zhong [(2003a), Proposition 3.2] for details.

We are ready to present the $\mathrm{P}_{\lambda_d}$-analogue of Lemma 5.2.

LEMMA 5.4. *For all $f : \mathbf{R}^d \to \mathbf{R}_+$ in $\mathcal{L}^1(\mathbf{R}^d) \cap \mathcal{L}^2(\mathbf{R}^d)$, and for all $\mu \in \mathcal{P}(\mathbf{R}_+)$,*

$$(5.13) \qquad 2^{-p} \|\widehat{f}\|^2_{\mathcal{L}^2(\mathrm{Q}_\mu^{\alpha p})} \leq \mathrm{E}_{\lambda_d}[|\mathrm{O}_\mu(f)|^2] \leq 2^{p\alpha} \|\widehat{f}\|^2_{\mathcal{L}^2(\mathrm{Q}_\mu^{\alpha p})},$$

*where the measure $\mathrm{Q}_\mu^\beta$ is defined in (5.10).*

PROOF. In the notation of the present note, if we further assume that $\widehat{f} \in \mathcal{L}^1(\mathbf{R}^d)$, then Lemma 3.5 of Khoshnevisan, Xiao and Zhong (2003a) and symmetry together show that

$$(5.14) \qquad \mathrm{E}_{\lambda_d}[|\mathrm{O}_\mu(f)|^2] = (2\pi)^{-d} \int_{\mathbf{R}^d} |\widehat{f}(\xi)|^2 \mathrm{E}[|\widehat{\mathrm{O}_\mu}(\xi)|^2] \, d\xi.$$

The lemma—under the extra assumption that $\widehat{f} \in \mathcal{L}^1(\mathbf{R}^d)$—follows from this, used in conjunction with (5.7) and (5.8). To drop the integrability condition on $\widehat{f}$, note a mollification argument reveals that all that is needed is $\widehat{f} \in \mathcal{L}^2(\mathbf{R}^d)$; but by the Plancherel theorem, this is equivalent to $f \in \mathcal{L}^2(\mathbf{R}^d)$. □

We are ready to dispense with the first part of the proof of Theorem 2.2.

PROOF OF THEOREM 2.2 (First half). Choose $\alpha \in (0,2]$ and an integer $p \geq 1$ such that

$$\alpha p = d - \beta. \tag{5.15}$$

Then we introduce an independent $p$-parameter additive $\alpha$-stable process $\vec{X} := \{\vec{X}(t)\}_{t \in \mathbf{R}_+^p}$ by

$$\vec{X}(t) := X_1(t_1) + \cdots + X_p(t_p) \qquad \forall t \in \mathbf{R}_+^p. \tag{5.16}$$

This also defines a $(1+p)$-parameter additive Lévy process $A := \{A(\mathbf{t})\}_{\mathbf{t} \in \mathbf{R}_+^{1+p}}$ defined by (5.2).

Now suppose there exists a $\mu \in \mathcal{P}(G)$ such that $\int_{\mathbf{R}^d} \mathcal{E}_{\chi_\xi}(\mu) \|\xi\|^{\beta - d} d\xi < +\infty$. Then, Lemma 5.2 and Plancherel's theorem, used in conjunction, tell us that there exists a (measurable) process $\{\ell_\mu(x)\}_{x \in \mathbf{R}^d}$ such that:

1. $\mathrm{E}[\|\ell_\mu\|_{\mathcal{L}^2(\mathbf{R}^d)}^2] = (2\pi)^{-d} \mathrm{E}[\|\widehat{\mathrm{O}_\mu}\|_{\mathcal{L}^2(\mathbf{R}^d)}^2] \leq 2^{p\alpha} \mathrm{Q}_\mu^{\alpha p}(\mathbf{R}^d) < +\infty$; see also Remark 5.3.
2. With probability one, for all bounded measurable functions $f : \mathbf{R}^d \to \mathbf{R}$, $\mathrm{O}_\mu(f) = \int_{\mathbf{R}^d} f(x) \ell_\mu(x) \, dx$.

Apply part 2 with $f(x) := \mathbf{1}_{A(G \times \mathbf{R}_+^p)}(x)$, and apply the Cauchy–Schwarz inequality to deduce that almost surely,

$$\begin{aligned}
1 = \mathrm{O}_\mu(\mathbf{1}_{A(G \times \mathbf{R}_+^p)}) &= \int_{\mathbf{R}^d} \sqrt{\mathbf{1}_{A(G \times \mathbf{R}_+^p)}(x)} \, \ell_\mu(x) \, dx \\
&\leq \sqrt{\lambda_d(A(G \times \mathbf{R}_+^p))} \|\ell_\mu\|_{\mathcal{L}^2(\mathbf{R}^d)}.
\end{aligned} \tag{5.17}$$

By the Cauchy–Schwarz inequality and part 1,

$$\mathrm{E}[\lambda_d(A(G \times \mathbf{R}_+^p))] \geq \frac{1}{\mathrm{E}[\|\ell_\mu\|_{\mathcal{L}^2(\mathbf{R}^d)}^2]} \geq \frac{1}{2^{p\alpha} \mathrm{Q}_\mu^{\alpha p}(\mathbf{R}^d)}. \tag{5.18}$$

Since $\mu \in \mathcal{P}(G)$ can be chosen arbitrarily as long as $\mathrm{Q}_\mu^{\alpha p}(\mathbf{R}^d) < +\infty$, and because of Remark 5.3 and (5.15), we have demonstrated that

$$\inf_{\mu \in \mathcal{P}(G)} \int_{\mathbf{R}^d} \mathcal{E}_{\chi_\xi}(\mu) \|\xi\|^{\beta - d} d\xi < +\infty \quad \Longrightarrow \quad \mathrm{E}[\lambda_d(A(G \times \mathbf{R}_+^p))] > 0. \tag{5.19}$$

According to Theorem 4.4, and thanks to (5.15), for any Borel set $F \subset \mathbf{R}^d$,

$$\mathrm{E}[\lambda_d(F \oplus \vec{X}(\mathbf{R}_+^p))] > 0 \quad \Longleftrightarrow \quad \mathcal{C}_\beta(F) > 0. \tag{5.20}$$

Since $X$ is independent of $\vec{X}$, we can apply this, conditionally, with $F := X(G)$, and then integrate $[d\mathrm{P}]$, to deduce that

$$\mathrm{E}[\lambda_d(A(G \times \mathbf{R}_+^p))] > 0 \quad \Longleftrightarrow \quad \mathrm{E}[\mathcal{C}_\beta(X(G))] > 0. \tag{5.21}$$

This and (5.19) together imply that

$$(5.22) \quad \inf_{\mu \in \mathcal{P}(G)} \int_{\mathbf{R}^d} \mathcal{E}_{\chi_\xi}(\mu) \|\xi\|^{\beta-d} \, d\xi < +\infty \quad \Longrightarrow \quad \mathrm{E}[\mathcal{C}_\beta(X(G))] > 0.$$

This proves fully half of Theorem 2.2. □

PROOF OF THEOREM 2.2 (Second half). We now prove the more difficult second half of Theorem 2.2; that is,

$$(5.23) \quad \mathrm{E}[\mathcal{C}_\beta(X(G))] > 0 \quad \Longrightarrow \quad \inf_{\mu \in \mathcal{P}(G)} \int_{\mathbf{R}^d} \mathcal{E}_{\chi_\xi}(\mu) \|\xi\|^{\beta-d} \, d\xi < +\infty.$$

In so doing, we can assume without loss of generality that the set $G$ is compact. Indeed, consider both sides, in (5.23), of "⇒" as set functions in $G$. Both of the said functions are Choquet capacities. Hence, Choquet's theorem reduces our analysis to the study of compact sets $G$.

Henceforth, $\{\varphi_\varepsilon\}_{\varepsilon>0}$ denotes the Gaussian approximation to the identity,

$$(5.24) \quad \varphi_\varepsilon(x) := (2\pi\varepsilon^2)^{-d/2} \exp\left(-\frac{\|x\|^2}{2\varepsilon^2}\right) \quad \forall x \in \mathbf{R}^d, \varepsilon > 0.$$

As we did earlier, we choose $\alpha \in (0, 2]$ and an integer $p \geq 1$ such that $\alpha p = d - \beta$. We bring in $p$ independent $\alpha$-stable Lévy processes $X_1, \ldots, X_p$, and construct the corresponding additive Lévy process $A := X \oplus \vec{X}$ defined by (5.2).

Let us start with setting some preliminary groundwork. To begin with, we define a $(1+p)$-parameter filtration $\mathfrak{F} := \{\mathfrak{F}(\mathbf{t})\}_{\mathbf{t} \in \mathbf{R}_+^{1+p}}$ by defining $\mathfrak{F}(\mathbf{t})$ to be the sigma-algebra defined by $\{A(\mathbf{r})\}_{\mathbf{r} \prec \mathbf{t}}$. Without loss of generality, we can assume that each $\mathfrak{F}(\mathbf{t})$ has been completed with respect to all measures $\mathrm{P}_x$ ($x \in \mathbf{R}^d$). We remark that $\mathfrak{F}$ is, indeed, a filtration in the partial order $\prec$. By this we mean that whenever $\mathbf{s} \prec \mathbf{t}$, then $\mathfrak{F}(\mathbf{s}) \subseteq \mathfrak{F}(\mathbf{t})$; a fact that can be readily checked.

Next we define, for any $\mu \in \mathcal{P}(G)$, the $(1+p)$-parameter process $\{M_\mu \varphi_\varepsilon(\mathbf{t})\}_{\mathbf{t} \in \mathbf{R}_+^{1+p}}$ as follows:

$$(5.25) \quad M_\mu \varphi_\varepsilon(\mathbf{t}) := \mathrm{E}_{\lambda_d}[\mathrm{O}_\mu(\varphi_\varepsilon) | \mathfrak{F}(\mathbf{t})] \quad \forall \mathbf{t} \in \mathbf{R}_+^{1+p},$$

where $\mathrm{O}_\mu(\varphi_\varepsilon)$ is defined by (5.3). It should be recognized that $M_\mu \varphi_\varepsilon$ is a $(1+p)$-parameter martingale in the partial order $\prec$ and in the infinite-measure space $(\Omega, \mathfrak{F}, \mathrm{P}_{\lambda_d})$. By a martingale here, we mean that whenever $\mathbf{s} \prec \mathbf{t}$, then $\mathrm{P}_{\lambda_d}$-almost surely,

$$(5.26) \quad \mathrm{E}_{\lambda_d}[M_\mu \varphi_\varepsilon(\mathbf{t}) | \mathfrak{F}(\mathbf{s})] = M_\mu \varphi_\varepsilon(\mathbf{s}).$$

By specializing Lemma 4.1 of Khoshnevisan, Xiao and Zhong (2003a) to the present setting, we obtain the following:

$$\mathrm{E}_{\lambda_d}[M_\mu \varphi_\varepsilon(\mathbf{t})] = 1 \qquad \forall\, \mathbf{t} \in \mathbf{R}_+^{1+p},$$

(5.27) $$\sup_{\mathbf{t} \in \mathbf{R}_+^{1+p}} \mathrm{E}_{\lambda_d}[(M_\mu \varphi_\varepsilon(\mathbf{t}))^2] \leq (2\pi)^{-d} \int_{\mathbf{R}^d} |\widehat{\varphi_\varepsilon}(\xi)|^2 \mathrm{E}[|\widehat{O_\mu}(\xi)|^2]\, d\xi$$

$$\leq 2^{p\alpha} \|\widehat{\varphi_\varepsilon}\|^2_{\mathcal{L}^2(Q_\mu^{\alpha p})};$$

see Lemma 5.1 for the last line.

Next, we work toward a bound in the reverse direction. For this, we note that for any $\mathbf{s} \in \mathbf{R}_+^{1+p}$,

(5.28) $$M_\mu \varphi_\varepsilon(\mathbf{s}) \geq \mathrm{E}_{\lambda_d}\bigg[\int_{\mathbf{t}\succ\mathbf{s}} \varphi_\varepsilon(A(\mathbf{t}))\boldsymbol{\mu}(d\mathbf{t})\bigg|\mathfrak{F}(\mathbf{s})\bigg] = \int_{\mathbf{t}\succ\mathbf{s}} P_{\mathbf{t}-\mathbf{s}}\varphi_\varepsilon(A(\mathbf{s}))\boldsymbol{\mu}(d\mathbf{t}),$$

where

(5.29) $$P_{\mathbf{t}}g(x) := \mathrm{E}[g(x + A(\mathbf{t}))] \qquad \forall\, \mathbf{t} \in \mathbf{R}_+^{1+p}, x \in \mathbf{R}^d,$$

and the last equality in (5.28) follows from the Markov property of the additive Lévy process $A$ under $\mathrm{P}_{\lambda_d}$. See Khoshnevisan and Xiao [(2002), Proposition 5.8] or Khoshnevisan, Xiao and Zhong [(2003a), Proposition 3.2].

Now suppose that $G \subset (0,\infty)$ is compact, and $\mathrm{E}[\mathcal{C}_\beta(X(G))] > 0$. By (5.21), this is equivalent to assuming

(5.30) $$\mathrm{E}[\lambda_d(A(G \times \mathbf{R}_+^p))] > 0.$$

By (5.28), $\mathrm{P}_{\lambda_d}$-almost surely,

(5.31) $$M_\mu \varphi_\varepsilon(\mathbf{s}) \geq \int_{\mathbf{t}\succ\mathbf{s}} P_{\mathbf{t}-\mathbf{s}}\varphi_\varepsilon(A(\mathbf{s}))\boldsymbol{\mu}(d\mathbf{t}) \cdot \mathbf{1}_{\{\|A(\mathbf{s})\|\leq\delta\}}$$

$$= \int_{\mathbf{t}\succ\mathbf{s}} \mathrm{E}[\varphi_\varepsilon(A(\mathbf{s}) + A'(\mathbf{t}-\mathbf{s}))|A(\mathbf{s})]\boldsymbol{\mu}(d\mathbf{t}) \cdot \mathbf{1}_{\{\|A(\mathbf{s})\|\leq\delta\}},$$

where $\{A'(\mathbf{t})\}_{\mathbf{t}\in\mathbf{R}_+^{1+p}}$ is an independent copy of $\{A(\mathbf{t})\}_{\mathbf{t}\in\mathbf{R}_+^{1+p}}$. In particular, $\mathrm{P}_{\lambda_d}$-almost surely,

(5.32) $$M_\mu \varphi_\varepsilon(\mathbf{s}) \geq \int_{\mathbf{t}\succ\mathbf{s}} \inf_{z\in\mathbf{R}^d:\|z\|\leq\delta} \mathrm{E}[\varphi_\varepsilon(z + A'(\mathbf{t}-\mathbf{s}))]\boldsymbol{\mu}(d\mathbf{t}) \cdot \mathbf{1}_{\{\|A(\mathbf{s})\|\leq\delta\}}.$$

On the other hand, one can directly check that

(5.33) $$\int_{\mathbf{t}\succ\mathbf{s}} \inf_{z\in\mathbf{R}^d:\|z\|\leq\delta} \mathrm{E}[\varphi_\varepsilon(z + A'(\mathbf{t}-\mathbf{s}))]\boldsymbol{\mu}(d\mathbf{t})$$

$$= \int_{s_0}^\infty \int_{s_1}^\infty \cdots \int_{s_p}^\infty \inf_{z\in\mathbf{R}^d:\|z\|\leq\delta} P_{\mathbf{t}-\mathbf{s}}\varphi_\varepsilon(z) e^{-\sum_{j=1}^p t_j} \mu(dt_0)\, d\vec{t}.$$

According to Lemma 3.1 of Khoshnevisan, Xiao and Zhong (2003a),

(5.34)
$$P_{\mathbf{t}-\mathbf{s}}\varphi_\varepsilon(0)$$
$$= (2\pi)^{-d} \int_{\mathbf{R}^d} e^{-(t_0-s_0)\Psi(\xi) - \sum_{j=1}^p (t_j - s_j)\|\xi\|^\alpha - (1/2)\varepsilon^2 \|\xi\|^2} \, d\xi.$$

Because the left-hand side is strictly positive, so is the right-hand side. In addition,

(5.35)
$$\int_{s_0}^\infty \int_{s_1}^\infty \cdots \int_{s_p}^\infty P_{\mathbf{t}-\mathbf{s}}\varphi_\varepsilon(0) e^{-\sum_{j=1}^p t_j} \mu(dt_0) \, d\vec{t}$$
$$= (2\pi)^{-d} \int_{s_0}^\infty \int_{\mathbf{R}^d} e^{-(1/2)\varepsilon^2 \|\xi\|^2 - (t_0 - s_0)\Psi(\xi)}$$
$$\times \prod_{j=1}^p \int_{s_j}^\infty e^{-(v - s_j)\|\xi\|^\alpha - v} \, dv \, d\xi \, \mu(dt_0)$$
$$= (2\pi)^{-d} e^{-\sum_{j=1}^p s_j} \int_{s_0}^\infty \int_{\mathbf{R}^d} e^{-(1/2)\varepsilon^2 \|\xi\|^2 - (t_0 - s_0)\Psi(\xi)} \frac{d\xi}{(1 + \|\xi\|^\alpha)^p} \mu(dt_0).$$

We plug this into (5.32) and deduce the following from Fatou's lemma: $P_{\lambda_d}$-almost surely, for all $\mathbf{s} \in \mathbb{Q}_+^{1+p}$,

(5.36)
$$M_\mu \varphi_\varepsilon(\mathbf{s}) \geq \mathbf{1}_{\{\|A(\mathbf{s})\| \leq \delta\}} (2\pi)^{-d} (1 + o(1)) e^{-\sum_{j=1}^p s_j}$$
$$\times \int_{s_0}^\infty \int_{\mathbf{R}^d} e^{-(1/2)\varepsilon^2 \|\xi\|^2 - (t_0 - s_0)\Psi(\xi)} \frac{d\xi}{(1 + \|\xi\|^\alpha)^p} \mu(dt_0),$$

where $o(1)$ is a term that goes to 0, uniformly in $\vec{s}$ and $\mu$ (but not $\varepsilon$), as $\delta \to 0$. (This follows merely from the Lipschitz continuity of $\varphi_\varepsilon$.)

For any $\delta > 0$, define $G^\delta$ to be the closed $\delta$-enlargement of $G$, and note that $G^\delta$ is compact in $\mathbf{R}_+$. Choose some point $\Delta \notin \mathbf{R}_+$, and let $T^{\delta,l}$ denote any measurable $(\mathbb{Q}_+ \cap G^\delta) \cup \Delta$-valued function on $\Omega$ such that $T^{\delta,l} \neq \Delta$ if and only if there exists some $\vec{t} \in [0,l]^p$ such that $\|A(T^{\delta,l}, \vec{t})\| \leq \delta$. This can always be done since the $X_j$'s have cadlag paths, and since $B(0, \delta) := \{x \in \mathbf{R}^d : \|x\| \leq \delta\}$ has an open interior. It may help to think, informally, that $T^{\delta,l}$ is any measurably selected point in $G^\delta$ such that, for some $\vec{t} \in [0,l]^p$, $\|A(T^{\delta,l}, \vec{t})\| \leq \delta$, as long as such a point exists. If such a point does not exist, then the value of $T^{\delta,l}$ is set to $\Delta$. [*Warning*: This is very close to, but *not* the same as, the construction of Khoshnevisan, Xiao and Zhong (2003a).] Thus, (5.36) implies that $P_{\lambda_d}$-almost surely,

(5.37)
$$\sup_{\mathbf{s} \in \mathbf{R}_+^{1+p}} M_\mu \varphi_\varepsilon(\mathbf{s}) \geq \mathbf{1}_{\{T^{\delta,l} \neq \Delta\}} \frac{(2\pi)^{-d}(1 + o(1))}{e^{pl}}$$
$$\times \int_{T^{\delta,l}}^\infty \int_{\mathbf{R}^d} \frac{e^{-(1/2)\varepsilon^2 \|\xi\|^2 - (t_0 - T^{\delta,l})\Psi(\xi)}}{(1 + \|\xi\|^\alpha)^p} \, d\xi \, \mu(dt_0).$$

Finally, we choose $\mu \in P(G^\delta)$ judiciously. Fix $l > 0$, and define

$$(5.38) \qquad \mu^{\delta,k}(\bullet) := \frac{P_{\lambda_d}\{T^{\delta,l} \in \bullet, T^{\delta,l} \neq \Delta, \|A(0)\| \leq k\}}{P_{\lambda_d}\{T^{\delta,l} \neq \Delta, \|A(0)\| \leq k\}}.$$

Then, thanks to (5.30), for all sufficiently large $k$, $\mu^{\delta,k} \in \mathcal{P}(G^\delta)$; see Khoshnevisan, Xiao and Zhong [(2003a), (4.3)] and its subsequent display. Furthermore, $P_{\lambda_d}$-almost surely,

$$(5.39) \quad \sup_{\mathbf{s} \in \mathbf{R}_+^{1+p}} M_{\mu^{\delta,k}} \varphi_\varepsilon(\mathbf{s}) \geq \mathbf{1}_{\{T^{\delta,l} \neq \Delta, \|A(0)\| \leq k\}} \frac{(2\pi)^{-d}(1 + o(1))}{e^{pl}}$$
$$\times \int_{T^{\delta,l}}^\infty \int_{\mathbf{R}^d} \frac{e^{-(1/2)\varepsilon^2 \|\xi\|^2 - (t_0 - T^{\delta,l})\Psi(\xi)}}{(1 + \|\xi\|^\alpha)^p} d\xi \, \mu^{\delta,k}(dt_0).$$

We can square both sides of this inequality, and then take expectations to deduce that

$$E_{\lambda_d}\left[\left(\sup_{\mathbf{s} \in \mathbf{R}_+^{1+p}} M_{\mu^{\delta,k}} \varphi_\varepsilon(\mathbf{s})\right)^2\right]$$
$$\geq P_{\lambda_d}\{T^{\delta,l} \neq \Delta, \|A(0)\| \leq k\} \times (2\pi)^{-2d}(1 + o(1))e^{-2pl}$$
$$(5.40) \qquad \times E_{\lambda_d}\left[\left(\int_{T^{\delta,l}}^\infty \int_{\mathbf{R}^d} \cdots d\xi \, \mu^{\delta,k}(dt_0)\right)^2 \Big| T^{\delta,l} \neq \Delta, \|A(0)\| \leq k\right]$$
$$= P_{\lambda_d}\{T^{\delta,l} \neq \Delta, \|A(0)\| \leq k\} \times (2\pi)^{-2d}(1 + o(1))e^{-2pl}$$
$$\times \int_0^\infty \left(\int_y^\infty \int_{\mathbf{R}^d} \frac{e^{-(1/2)\varepsilon^2 \|\xi\|^2 - (x-y)\Psi(\xi)}}{(1 + \|\xi\|^\alpha)^p} d\xi \, \mu^{\delta,k}(dx)\right)^2 \mu^{\delta,k}(dy).$$

From the Cauchy–Schwarz inequality, after making an appeal to the fact that in the integrand $x \geq y$, we can deduce the following:

$$E_{\lambda_d}\left[\left(\sup_{\mathbf{s} \in \mathbf{R}_+^{1+p}} M_{\mu^{\delta,k}} \varphi_\varepsilon(\mathbf{s})\right)^2\right]$$
$$(5.41) \qquad \geq P_{\lambda_d}\{T^{\delta,l} \neq \Delta, \|A(0)\| \leq k\} \times (2\pi)^{-2d}(1 + o(1))e^{-2pl}$$
$$\times \left(\int_0^\infty \int_y^\infty \int_{\mathbf{R}^d} e^{-(1/2)\varepsilon^2 \|\xi\|^2 - |x-y|\Psi(\mathrm{sgn}(x-y)\xi)}\right.$$
$$\left. \times [(1 + \|\xi\|^\alpha)^p]^{-1} d\xi \, \mu^{\delta,k}(dx) \mu^{\delta,k}(dy)\right)^2.$$

This time, $o(1)$ is a term that goes to 0, uniformly over all $k \geq 1$, as $\delta \to 0$. We intend to show that, in the preceding display, we can replace, at little cost, $\int_y^\infty$ by $\int_0^\infty$. In order to do this, we need some preliminary setup. Most significantly, we need a new partial order on the enlarged parameter space $\mathbf{R}_+^{1+p}$.

For any $\mathbf{s}, \mathbf{t} \in \mathbf{R}_+^{1+p}$, we write

$$(5.42) \qquad \mathbf{s} \preccurlyeq \mathbf{t} \iff s_0 \geq t_0 \quad \text{but for all } j = 1, \ldots, p,\ s_j \leq t_j.$$

This is an entirely different partial order from $\prec$, and gives rise to a new $(1+p)$-parameter filtration $\mathfrak{R} := \{\mathfrak{R}(\mathbf{t})\}_{\mathbf{t} \in \mathbf{R}_+^{1+p}}$, where $\mathfrak{R}(\mathbf{t})$ is defined to be the sigma-algebra generated by $\{A(\mathbf{r})\}_{\mathbf{r} \preccurlyeq \mathbf{t}}$. Without loss of generality, we can assume that each $\mathfrak{R}(\mathbf{t})$ is complete with respect to every $\mathrm{P}_x$ ($x \in \mathbf{R}^d$). As we did for $\mathfrak{F}$, we remark that $\mathfrak{R}$ is a filtration in the new partial order $\preccurlyeq$ and, under the $\sigma$-finite measure $\mathrm{P}_{\lambda_d}$, $X$ satisfies the Markov property with respect to $\mathfrak{R}$. (The Fraktur letters $\mathfrak{F}$ and $\mathfrak{R}$ are chosen to remind the reader of "forward" and "reverse," since they refer to the time-order of the process $X$.)

Consider the $(1+p)$-parameter process $N_{\mu^{k,\delta}} \varphi_\varepsilon := \{N_{\mu^{k,\delta}} \varphi_\varepsilon(\mathbf{t})\}_{\mathbf{t} \in \mathbf{R}_+^{1+p}}$ that is defined by the following:

$$(5.43) \qquad N_{\mu^{k,\delta}} \varphi_\varepsilon(\mathbf{t}) := \mathrm{E}_{\lambda_d}[\mathrm{O}_{\mu^{k,\delta}}(\varphi_\varepsilon) | \mathfrak{R}(\mathbf{t})].$$

Clearly, this is a martingale in the partial order $\preccurlyeq$.

By using a similar argument as that which led to (5.41), we arrive at the following (here, it is essential to work with the infinite measure $\mathrm{P}_{\lambda_d}$ instead of P):

$$
\begin{aligned}
&\mathrm{E}_{\lambda_d}\left[\left(\sup_{\mathbf{s} \in \mathbf{R}_+^{1+p}} N_{\mu^{\delta,k}} \varphi_\varepsilon(\mathbf{s})\right)^2\right] \\
&\qquad \geq \mathrm{P}_{\lambda_d}\{T^{\delta,l} \neq \Delta, \|A(0)\| \leq k\} \times (2\pi)^{-2d}(1+o(1))e^{-2pl} \\
&\qquad \times \left(\int_0^\infty \int_0^y \int_{\mathbf{R}^d} e^{-(1/2)\varepsilon^2 \|\xi\|^2 - |x-y|\Psi(\mathrm{sgn}(x-y)\xi)} \right. \\
&\qquad\qquad \left. \times [(1+\|\xi\|^\alpha)^p]^{-1} d\xi\, \mu^{\delta,k}(dx) \mu^{\delta,k}(dy)\right)^2.
\end{aligned}
$$
(5.44)

[En route, this shows that the terms inside $(\cdots)^2$ are nonnegative real too.] Thus, we add (5.41) and (5.44), and use $2(a^2 + b^2) \geq (a+b)^2$—valid for all real $a, b$—to obtain the following:

$$
\begin{aligned}
&\mathrm{E}_{\lambda_d}\left[\left(\sup_{\mathbf{s} \in \mathbf{R}_+^{1+p}} M_{\mu^{\delta,k}} \varphi_\varepsilon(\mathbf{s})\right)^2\right] + \mathrm{E}_{\lambda_d}\left[\left(\sup_{\mathbf{s} \in \mathbf{R}_+^{1+p}} N_{\mu^{\delta,k}} \varphi_\varepsilon(\mathbf{s})\right)^2\right] \\
&\qquad \geq \mathrm{P}_{\lambda_d}\{T^{\delta,l} \neq \Delta, \|A(0)\| \leq k\} \times \left(\frac{1}{2(2\pi)^{2d}} + o(1)\right) e^{-2pl}
\end{aligned}
$$
(5.45)

$$\times \left( \int_0^\infty \int_0^\infty \int_{\mathbf{R}^d} e^{-(1/2)\varepsilon^2 \|\xi\|^2 - |x-y|\Psi(\operatorname{sgn}(x-y)\xi)} \right.$$
$$\left. \times [(1 + \|\xi\|^\alpha)^p]^{-1} d\xi \, \mu^{\delta,k}(dx) \mu^{\delta,k}(dy) \right)^2.$$

Now the integrand is absolutely integrable $[d\xi \times d\mu^{k,\delta} \times d\mu^{k,\delta}]$. Thus, by the Fubini–Tonelli theorem, we can interchange the order of the integrals, and obtain the following:

$$\int_0^\infty \int_0^\infty \int_{\mathbf{R}^d} \frac{e^{-(1/2)\varepsilon^2 \|\xi\|^2 - |x-y|\Psi(\operatorname{sgn}(x-y)\xi)}}{(1 + \|\xi\|^\alpha)^p} d\xi \, \mu^{\delta,k}(dx) \mu^{\delta,k}(dy)$$

(5.46)
$$= \int_{\mathbf{R}^d} \frac{e^{-(1/2)\varepsilon^2 \|\xi\|^2}}{(1 + \|\xi\|^\alpha)^p} \mathcal{E}_{\chi_\xi}(\mu^{\delta,k}) \, d\xi$$

$$\geq 2^{-\alpha p}(2\pi)^d \int_{\mathbf{R}^d} e^{-(1/2)\varepsilon^2 \|\xi\|^2} Q^{\alpha p}_{\mu^{k,\delta}}(d\xi) \geq 2^{-\alpha p}(2\pi)^d \|\widehat{\varphi_\varepsilon}\|^2_{\mathcal{L}^2(Q^{\alpha p}_{\mu^{k,\delta}})}.$$

In the above, the first inequality follows from (5.8) and (5.10), and the second inequality follows from $0 < \widehat{\varphi_\varepsilon}(\xi) \leq 1$.

In other words, after recalling (5.10), we arrive at the following:

$$\mathrm{E}_{\lambda_d}\left[\left(\sup_{\mathbf{s} \in \mathbf{R}^{1+p}_+} M_{\mu^{\delta,k}} \varphi_\varepsilon(\mathbf{s})\right)^2\right] + \mathrm{E}_{\lambda_d}\left[\left(\sup_{\mathbf{s} \in \mathbf{R}^{1+p}_+} N_{\mu^{\delta,k}} \varphi_\varepsilon(\mathbf{s})\right)^2\right]$$
(5.47)
$$\geq \mathrm{P}_{\lambda_d}\{T^{\delta,l} \neq \Delta, \|A(0)\| \leq k\} \times \frac{1 + o(1)}{e^{2pl} 2^{1+2\alpha p}} \|\widehat{\varphi_\varepsilon}\|^4_{\mathcal{L}^2(Q^{\alpha p}_{\mu^{\delta,k}})}.$$

We recall that $o(1)$ is a term that tends to zero, as $\delta \to 0$, uniformly in all of the variables except $\varepsilon$.

It turns out that, under the infinite measure $\mathrm{P}_{\lambda_d}$, both filtrations $\mathfrak{F}$ and $\mathfrak{R}$ are commuting in the sense of Khoshnevisan [(2002), page 233]; see Khoshnevisan, Xiao and Zhong [(2003a), proof of Lemma 4.2, page 1111] for a discussion of a much more general property. Thus, by the Cairoli inequality [Khoshnevisan (2002), Theorem 2.3.2, Chapter 7],

$$\mathrm{E}_{\lambda_d}\left[\left(\sup_{\mathbf{s} \in \mathbf{R}^{1+p}_+} M_{\mu^{\delta,k}} \varphi_\varepsilon(\mathbf{s})\right)^2\right] \leq 4^{p+1} \sup_{\mathbf{s} \in \mathbf{R}^{1+p}_+} \mathrm{E}_{\lambda_d}[M^2_{\mu^{\delta,k}}(\mathbf{s})]$$
(5.48)
$$= 4^{p+1} \sup_{\mathbf{s} \in \mathbf{R}^{1+p}_+} \mathrm{E}_{\lambda_d}[(O_{\mu^{\delta,k}}(\varphi_\varepsilon))^2]$$

$$\leq 4^{p+1} 2^{p\alpha} \|\widehat{\varphi_\varepsilon}\|^2_{\mathcal{L}^2(Q^{\alpha p}_{\mu^{\delta,k}})}$$

[cf. (5.25) and (5.27); the fact that $\mathrm{P}_{\lambda_d}$ is not a probability measure does not cause any difficulties here]. Moreover, the preceding remains valid if we replace $M_{\mu^{k,\delta}}$ by $N_{\mu^{k,\delta}}$ everywhere. Thus, solving the preceding two displays leads us to the following:

$$(5.49) \quad P_{\lambda_d}\{T^{\delta,l} \neq \Delta, \|A(0)\| \leq k\} \leq e^{2pl} 2^{3(1+\alpha p)+2p}(1+o(1))\|\widehat{\varphi_\varepsilon}\|^2_{\mathcal{L}^2(Q^{\alpha p}_{\mu^{\delta,k}})}.$$

Now we let $k \to \infty$ and appeal to Fatou's lemma to see that there must exist $\mu^\delta \in \mathcal{P}(G^\delta)$ such that

$$(5.50) \quad P_{\lambda_d}\{T^{\delta,l} \neq \Delta\} \leq e^{2pl} 2^{3(1+\alpha p)+2p}(1+o(1))\|\widehat{\varphi_\varepsilon}\|^{-2}_{\mathcal{L}^2(Q^{\alpha p}_{\mu^\delta})}.$$

In order to deduce the above, note that all of the probability measures $\{\mu^{k,\delta}\}_{k \geq 1}$ live on the same compact set $G^\delta$. Therefore, we can extract a subsequence that converges weakly to $\mu^\delta \in \mathcal{P}(G^\delta)$. To finish, note that $|\widehat{\varphi_\varepsilon}(\xi)|^2 = \exp(-\varepsilon^2\|\xi\|^2)$ is a bounded continuous function of $\xi$ and it is in $L^1(\mathbf{R}^d)$. Hence, by the Fubini–Tonelli theorem, we have

$$(5.51) \quad \|\widehat{\varphi_\varepsilon}\|^2_{\mathcal{L}^2(Q^{\alpha p}_{\mu^{\delta,k}})} = \iint \left[\int_{\mathbf{R}^d} e^{-\varepsilon^2\|\xi\|^2} \frac{e^{-|s-t|\Psi(\operatorname{sgn}(s-t)\xi)}}{(1+\|\xi\|)^{\alpha p}} d\xi\right] \mu^{\delta,k}(ds)\mu^{\delta,k}(dt),$$

and the kernel in the brackets is a bounded continuous function of $(s,t)$. So we obtain the asserted bound in (5.50).

Next, we let $\delta \downarrow 0$ in (5.50), and appeal to Fatou's lemma and compactness once more in order to obtain the following: There exists $\mu \in \mathcal{P}(G)$ such that

$$(5.52) \quad P_{\lambda_d}\{0 \in \overline{A(G \times [0,l]^p)}\} \leq c_{p,l,\alpha}\|\widehat{\varphi_\varepsilon}\|^{-2}_{\mathcal{L}^2(Q^{\alpha p}_\mu)},$$

where $c_{p,l,\alpha} := e^{2pl} 2^{3(1+\alpha p)+2p}$. We can now let $\varepsilon \downarrow 0$, and appeal to the monotone convergence theorem, to see that

$$(5.53) \quad \begin{aligned} E[\lambda_d(A(G \times [0,l]^p))] &\leq P_{\lambda_d}\{0 \in \overline{A(G \times [0,l]^p)}\} \\ &\leq c_{p,l,\alpha} \lim_{\varepsilon \to 0} \|\widehat{\varphi_\varepsilon}\|^{-2}_{\mathcal{L}^2(Q^{\alpha p}_\mu)} = \frac{c_{p,l,\alpha}}{Q^{\alpha p}_\mu(\mathbf{R}^d)}.\end{aligned}$$

In accordance with Remark 5.3, for this choice of $\mu \in \mathcal{P}(G)$, we have the following:

$$(5.54) \quad E[\lambda_d(A(G \times [0,l]^p))] > 0 \iff \int_{\mathbf{R}^d} \mathcal{E}_{\chi_\xi}(\mu)\|\xi\|^{-\alpha p} d\xi < +\infty.$$

At this stage, we can apply Theorem 4.4 [see also (5.20)], conditionally, with $F := X(G)$ to deduce that (5.30) holds if and only if $\mathcal{C}_{d-\alpha p}(F) > 0$ with positive probability. Hence, we have proven (5.23), and this completes our proof. $\square$

REMARK 5.5. Our proof of Theorem 2.2 is a self-contained argument for deriving the following:

$$(5.55) \quad \mathrm{E}[\mathcal{C}_\beta(X(G))] > 0 \iff \inf_{\mu \in \mathcal{P}(G)} \int_{\mathbf{R}^d} \mathcal{E}_{\chi_\xi}(\mu) \|\xi\|^{\beta-d}\, d\xi < +\infty.$$

We can then use Proposition 2.8 to conclude that the preceding is also equivalent to the condition that $X(G)$ almost surely has positive $\beta$-dimensional Bessel–Riesz capacity.

In this remark, we describe a proof that (5.55) implies Proposition 2.8. To do so, we need only to prove that

$$(5.56) \quad \inf_{\mu \in \mathcal{P}(G)} \int_{\mathbf{R}^d} \mathcal{E}_{\chi_\xi}(\mu) \|\xi\|^{\beta-d}\, d\xi < +\infty \implies \mathcal{C}_\beta(X(G)) > 0 \quad \text{a.s.}$$

That is, we assume that there exists $\mu \in \mathcal{P}(G)$ such that $\int_{\mathbf{R}^d} \mathcal{E}_{\chi_\xi}(\mu) \|\xi\|^{\beta-d}\, d\xi < +\infty$, and prove that, with probability one, $\mathcal{C}_\beta(X(G)) > 0$.

For such a probability measure $\mu \in \mathcal{P}(G)$, we define the occupation measure $\Lambda_\mu$ of $X$ by $\Lambda_\mu(A) := \int \mathbf{1}_A(X(s))\mu(ds)$, for all Borel sets $A \subset \mathbf{R}^d$. Informally, this is *exactly* the same as $\mathrm{O}_\mu(\mathbf{1}_A)$, where $p = 0$; see (5.3). Note that $\Lambda_\mu \in \mathcal{P}(X(G))$ a.s. and thanks to Plancherel's theorem in the form of (7.22) below, there exists a constant $c'_{d,\beta}$ such that $\mathcal{E}_\beta(\Lambda_\mu) = c'_{d,\beta} \int_{\mathbf{R}^d} |\widehat{\Lambda_\mu}(\xi)|^2 \|\xi\|^{\beta-d}\, d\xi$. On the other hand, by (5.7) (with $p := 0$), $\mathrm{E}[|\widehat{\Lambda_\mu}(\xi)|^2] = \mathcal{E}_{\chi_\xi}(\mu)$. Thus,

$$(5.57) \quad \mathrm{E}[\mathcal{E}_\beta(\Lambda_\mu)] = c'_{d,\beta} \int_{\mathbf{R}^d} \mathcal{E}_{\chi_\xi}(\mu) \|\xi\|^{\beta-d}\, d\xi,$$

which is finite. Thus, $\mathcal{E}_\beta(\Lambda_\mu)$ is finite almost surely, whence (5.56).

**6. Kahane's problem for self-intersections.** We now return to Kahane's problems, mentioned in the Introduction, regarding when $X(F) \cap X(G) \neq \varnothing$ for disjoint sets $F$ and $G$ in $\mathbf{R}_+$. The following is the most general answer that we have been able to find.

THEOREM 6.1. *If $X$ is a Lévy process in $\mathbf{R}^d$ with Lévy exponent $\Psi$, then given any two disjoint Borel sets $F, G \subset \mathbf{R}_+$, $\mathrm{E}[\lambda_d(X(F) \ominus X(G))] > 0$ if and only if there exists $\mu \in \mathcal{P}(F \times G)$ such that*

$$(6.1) \quad \begin{aligned} & \int_{\mathbf{R}^d} \mathcal{E}_{\chi_\xi \otimes \chi_{-\xi}}(\mu)\, d\xi \\ & := \int_{\mathbf{R}^d} \iint \chi_\xi(s_1 - t_1)\chi_{-\xi}(s_2 - t_2)\mu(ds)\mu(dt)\, d\xi < +\infty. \end{aligned}$$

*If, in addition, the distribution of $X(t)$ is equivalent to $\lambda_d$ for all $t > 0$, then the above condition* (6.1) *is also equivalent to* $\mathrm{P}\{X(F) \cap X(G) \neq \varnothing\} > 0$.

In the symmetric case, $\chi_\xi$ is real and positive. So by the Fubini–Tonelli theorem, we have the following:

COROLLARY 6.2 (Kahane's problem). *Let $X$ be a symmetric Lévy process in $\mathbf{R}^d$ with Lévy exponent $\Psi$. If the distribution of $X(t)$ is equivalent to $\lambda_d$ for all $t > 0$, then $\mathrm{P}\{X(F) \cap X(G) \neq \varnothing\} > 0$ if and only if $\mathcal{C}_f(F \times G) > 0$, where for all $x \in \mathbf{R}^2$,*

$$(6.2) \qquad f(x) := \int_{\mathbf{R}^d} \chi_\xi \otimes \chi_\xi(x) \, d\xi := \int_{\mathbf{R}^d} e^{-(|x_1|+|x_2|)\Psi(\xi)} \, d\xi.$$

EXAMPLE 6.3. If $X$ is a symmetric $\alpha$-stable Lévy process in $\mathbf{R}^d$, then $\Psi(\xi) \geq 0$ and $c\|\xi\|^\alpha \leq \Psi(\xi) \leq C\|\xi\|^\alpha$ for some constants $0 < c \leq C$, and we readily obtain the following consequence which solves the problem, due to Kahane, mentioned in the Introduction:

$$(6.3) \qquad \mathrm{P}\{X(F) \cap X(G) \neq \varnothing\} > 0 \quad \Longleftrightarrow \quad \mathcal{C}_{d/\alpha}(F \times G) > 0.$$

This was previously known only when $\alpha = 2$; that is, when $X$ is a Brownian motion [Khoshnevisan (1999), Theorem 8.2].

Now we begin proving our way toward Theorem 6.1. The first step is a simplification that is well known, as well as interesting on its own. Namely, in order to prove Theorem 6.1, it suffices to prove the following:

THEOREM 6.4. *Suppose $X_1$ and $X_2$ are independent Lévy process in $\mathbf{R}^d$ with Lévy exponents $\Psi_1$ and $\Psi_2$. Then, given any two Borel sets $F_1$ and $F_2$, both in $\mathbf{R}_+$, $\mathrm{E}[\lambda_d(X_1(F_1) \ominus X_2(F_2))] > 0$ if and only if there exists $\mu \in \mathcal{P}(F_1 \times F_2)$ such that*

$$(6.4) \qquad \int_{\mathbf{R}^d} \mathcal{E}_{\chi_\xi^{\Psi_1} \otimes \chi_{-\xi}^{\Psi_2}}(\mu) \, d\xi < +\infty.$$

*If, in addition, the distribution of $X_1 \ominus X_2(\vec{t})$ is equivalent to $\lambda_d$ for all $\vec{t} \in (0, \infty)^2$, then the above condition is also equivalent to the condition that $\mathrm{P}\{X_1(F_1) \cap X_2(F_2) \neq \varnothing\} > 0$.*

PROOF. Consider the two-parameter additive Lévy process $\vec{X} := X_1 \ominus X_2$. The Lévy exponent of the process $\vec{X}$ is the function $\vec{\Psi}(\xi) := (\Psi_1(\xi), \Psi_2(-\xi))$; of course, $\Psi_2(-\xi)$ is the complex conjugate of $\Psi_2(\xi)$. The necessary and sufficient condition for the positivity of $\mathrm{E}[\lambda_d(X_1(F_1) \ominus X_2(F_2))]$ follows from Theorem 4.8. To finish, we can apply Lemma 4.1 with the choices, $F := \{0\}$ and $G := F_1 \times F_2$. □

We are finally ready to prove Theorem 6.1.

PROOF OF THEOREM 6.1. It suffices to prove this theorem for $F$ and $G$ compact subsets of $\mathbf{R}_+$.

We can simplify the problem further by assuming, without loss of generality, that there exist $0 < a < b < c < d$ such that $F \subset [a,b]$ and $G \subset [c,d]$. Choose any nonrandom number $\tau \in (b,c)$, and note that the translation invariance of $\lambda_d$ implies that $\lambda_d(X(F) \ominus X(G)) > 0$ if and only if $\lambda_d(X_1(F) \ominus X_2(G \ominus \tau)) > 0$, where $X_1(t) := X(t)$ $(0 \le t \le \tau)$ and $X_2(t) := X(t+\tau) - X(\tau)$ $(t \ge 0)$. Clearly, $X_1$ and $X_2$ are independent Lévy processes both when exponent $\Psi$ is verified automatically. Thus, by Theorem 6.4,

$$
\begin{aligned}
&\mathrm{E}[\lambda_d(X(F) \ominus X(G))] > 0 \\
&\quad \Longleftrightarrow \quad \inf_{\mu \in \mathcal{P}(F \times G \ominus \tau)} \int_{\mathbf{R}^d} \mathcal{E}_{\chi_\xi \otimes \chi_{-\xi}}(\mu)\, d\xi < +\infty.
\end{aligned}
\tag{6.5}
$$

By the explicit form of the latter energies [cf. (6.1)], the above condition is equivalent to the existence of $\nu \in \mathcal{P}(F \times G)$ such that $\int_{\mathbf{R}^d} \mathcal{E}_{\chi_\xi \otimes \chi_{-\xi}}(\nu)\, d\xi$ is finite. This proves the first half of the theorem.

Now suppose, in addition, that the distribution of $X(t)$ is equivalent to $\lambda_d$ for all $t > 0$. Consider the two-parameter additive Lévy process $\vec{X} := X_1 \ominus X_2$, where $X_1$ and $X_2$ are the same processes we used earlier in this proof, and note that the distribution of $\vec{X}(\vec{t})$ is equivalent to $\lambda_d$ for all $\vec{t} := (t_1, t_2) \in (0, \infty)^2$. Hence, Lemma 4.1 implies that

$$
\begin{aligned}
&\mathrm{P}\{X_1(F) \cap X_2(G \ominus \tau) \ne \varnothing\} > 0 \\
&\quad \Longleftrightarrow \quad \inf_{\mu \in \mathcal{P}(F \times G)} \int_{\mathbf{R}^d} \mathcal{E}_{\chi_\xi \otimes \chi_{-\xi}}(\mu)\, d\xi < +\infty.
\end{aligned}
\tag{6.6}
$$

Equivalently,

$$
\begin{aligned}
&\mathrm{P}\{X_1(F) \cap X_2(G \ominus \tau) \ne \varnothing | \mathfrak{X}(\tau)\} > 0 \quad \text{with positive probability} \\
&\quad \Longleftrightarrow \quad \inf_{\mu \in \mathcal{P}(F \times G)} \int_{\mathbf{R}^d} \mathcal{E}_{\chi_\xi \otimes \chi_{-\xi}}(\mu)\, d\xi < +\infty,
\end{aligned}
\tag{6.7}
$$

where $\mathfrak{X}(\tau)$ denotes the sigma-algebra generated by $\{X(u);\ u \in [0,\tau]\}$. Now, $X_2$ is independent of $\mathfrak{X}(\tau)$ and $X_1$. So we can apply Lemma 4.1 [with $p := 1$, $F$ replaced with $Z \oplus X_1(F)$, and $\vec{X}$ replaced with $X_2$] to deduce that for any a.s.-finite $\mathfrak{X}(\tau)$-measurable random variable $Z$,

$$
\begin{aligned}
&\mathrm{P}\{Z \oplus X_1(F) \cap X_2(G \ominus \tau) \ne \varnothing | \mathfrak{X}(\tau)\} > 0 \quad \text{with positive probability} \\
&\quad \Longleftrightarrow \quad \inf_{\mu \in \mathcal{P}(F \times G)} \int_{\mathbf{R}^d} \mathcal{E}_{\chi_\xi \otimes \chi_{-\xi}}(\mu)\, d\xi < +\infty.
\end{aligned}
\tag{6.8}
$$

Choose $Z := -X(\tau)$ and unscramble the above to conclude the proof. $\square$

Kahane (1983) has also studied the existence of self-intersections of a symmetric stable Lévy process $X = \{X(t)\}$ when $t$ is restricted to the disjoint compact subsets $E_1, E_2, \ldots, E_k$ of $\mathbf{R}_+$ ($k \geq 2$). The proof of Theorem 6.1 can be modified to give a necessary and sufficient condition for $P\{X(E_1) \cap \cdots \cap X(E_k) \neq \varnothing\} > 0$. For simplicity, we content ourselves by deriving the following result from Theorem 4.2 under the extra assumption that $X$ is symmetric and absolutely continuous. We point out that when $k = 2$, the conditions of Theorem 6.5 and Corollary 6.2 are not always comparable.

THEOREM 6.5. *Let $X$ be a symmetric Lévy process in $\mathbf{R}^d$ with Lévy exponent $\Psi$. Suppose that, for every fixed $t > 0$, $e^{-t\Psi(\cdot)} \in \mathcal{L}^1(\mathbf{R}^d)$. Then, for all disjoint compact sets $E_1, \ldots, E_k \subset \mathbf{R}_+$, $P\{X(E_1) \cap \cdots \cap X(E_k) \neq \varnothing\} > 0$ if and only if $\mathcal{C}_f(E_1 \times E_2 \times \cdots \times E_k) > 0$. Here,*

$$f(x) := (2\pi)^{-d(k-1)} \int_{\mathbf{R}^{d(k-1)}} \exp\left(-\sum_{j=1}^{k} |x_j| \Psi(\xi_{j-1} - \xi_j)\right) d\xi$$

(6.9)
$$\forall x \in \mathbf{R}^k.$$

*We have written $\xi \in \mathbf{R}^{d(k-1)}$ as $\xi := \xi_1 \otimes \cdots \otimes \xi_{k-1}$, where $\xi_j \in \mathbf{R}^d$. In addition, $\xi_0 := \xi_k := 0$.*

PROOF. By the proof of Theorem 6.1, it suffices to consider $k$ independent symmetric Lévy processes $X_1, \ldots, X_k$ in $\mathbf{R}^d$ with exponent $\Psi$. We define a multiparameter process $\vec{X} := \{\vec{X}(t)\}_{t \in \mathbf{R}_+^k}$, with values in $\mathbf{R}^{(k-1)d}$, by

(6.10) $\qquad \vec{X}(t) = (X_2(t_2) - X_1(t_1), \ldots, X_k(t_k) - X_{k-1}(t_{k-1}))$.

Then $\vec{X}$ can be expressed as an additive Lévy process in $\mathbf{R}^{(k-1)d}$ with Lévy exponent $(\Psi_1, \ldots, \Psi_k)$, where for every $j = 1, \ldots, k$, $\Psi_j$ is defined by

(6.11) $\qquad \Psi_j(\xi) = \Psi(\xi_{j-1} - \xi_j) \qquad \forall \xi = (\xi_1, \ldots, \xi_{k-1}) \in \mathbf{R}^{(k-1)d}$.

It is easy to verify that, under our assumptions, $\vec{X}$ is a symmetric and absolutely continuous additive Lévy process whose gauge function is given by (6.9). Because $P\{\bigcap_{j=1}^k X_j(E_j) \neq \varnothing\} > 0$ if and only if $P\{\vec{X}^{-1}(0) \cap (E_1 \times \cdots \times E_k) \neq \varnothing\} > 0$, Theorem 6.5 follows from Theorem 4.2 and Remark 4.3. Related information can be found in Khoshnevisan and Xiao [(2002), pages 93 and 94]. $\square$

**7. Examples of capacity and dimension computations.**

7.1. *Isotropic processes*: *image.* Throughout this section we consider an isotropic Lévy process $X := \{X(t)\}_{t \geq 0}$ with an exponent $\Psi$ that is regularly varying at infinity with index $\alpha \in (0, 2]$. Thus, we may write

$$\Psi(\xi) = \|\xi\|^\alpha \kappa(\|\xi\|) \qquad \forall \xi \in \mathbf{R}^d \setminus \{0\}. \tag{7.1}$$

Here, $\kappa : (0, \infty) \to \mathbf{R}_+$ is a function that is slowly varying at infinity. We now derive the following application of Theorem 2.2 for a broad class of such processes.

THEOREM 7.1. *Suppose $\kappa : (0, \infty) \to \mathbf{R}_+$ is continuous and slowly varying at infinity. Then, for any nonrandom Borel set $G \subset \mathbf{R}_+$, and all $\beta \in (0, d)$,*

$$\mathcal{C}_\beta(X(G)) = 0 \quad a.s. \iff \mathcal{C}_{g_\kappa}(G) = 0, \tag{7.2}$$

*where*

$$g_\kappa(x) := |x|^{-\beta/\alpha} [\kappa^\#(|x|^{-1/\alpha})]^\beta. \tag{7.3}$$

*Here, $\kappa^\#$ is the* de Bruijin conjugate *of $\kappa$.*

REMARK 7.2. It is known that $\kappa^\#$ is a slowly varying function [Bingham, Goldie and Teugels (1987), Theorem 1.5.13]. In many cases, the function $\kappa^\#$ can be estimated and/or computed with great accuracy; see Bingham, Goldie and Teugels [(1987), Section 5.2 and Appendix 5].

REMARK 7.3. If, in Theorem 7.1, we further assume that the function $\kappa(e^t)$ is regularly varying at infinity, then we can choose $g_\kappa$ as follows:

$$g_\kappa(x) := |x|^{-\beta/\alpha} \left[\kappa\left(\frac{1}{|x|}\right)\right]^{-\beta}. \tag{7.4}$$

The proof of this will be given in Remark 7.6 below.

Because $\mathcal{C}_f$ is determined by the behavior of $f$ at the origin, Theorem 7.1 follows from Corollary 2.4 at once if we could prove that

$$0 < \liminf_{|x| \to 0} \frac{f_{d-\beta}(x)}{g_\kappa(x)} \leq \limsup_{|x| \to 0} \frac{f_{d-\beta}(x)}{g_\kappa(x)} < +\infty. \tag{7.5}$$

Recall (2.7), integrate by parts, and change variables to see that

$$|x|^{\beta/\alpha} f_{d-\beta}(x) = v_d \int_0^\infty e^{-r^\alpha \kappa(r|x|^{-1/\alpha})} r^{\beta-1} \, dr, \tag{7.6}$$

where $v_d := \lambda_{d-1}(S^{d-1})$. Our next lemma will be used to describe the asymptotic behavior of $f_{d-\beta}(x)$ for $x$ near zero. We adopt the following notation: Given two nonnegative functions $h$ and $g$, and $x_0 \in [0, \infty]$, we write $h(x) \asymp g(x)$ $(x \to x_0)$ to mean that there exists a neighborhood $N$ of $x_0$ such that uniformly for $x \in N$, the ratio of $h(x)$ to $g(x)$ is bounded away from zero and infinity. (If $x_0 = +\infty$, neighborhood holds in the sense of the one-point compactification of $\mathbf{R}_+$. If no range of $x$ is specified, then the inequality holds for all $x$.)

LEMMA 7.4. *Under the conditions of Theorem 7.1, for any $\beta > 0$,*

$$(7.7) \qquad \int_0^\infty e^{-r^\alpha \kappa(nr)} r^{\beta-1}\, dr \asymp \varepsilon_n^\beta \qquad (n \to \infty),$$

*where $\varepsilon_n$ is any solution to $\varepsilon_n^\alpha \kappa(n\varepsilon_n) = 1$.*

PROOF. Let us begin by proving the existence of $\{\varepsilon_n\}_{n \geq 1}$. It follows from (7.1) that, for every fixed integer $n \geq 1$, $\lim_{x \to 0} x^\alpha \kappa(nx) = 0$ [since $\Psi(0) = 0$] and $\lim_{x \to \infty} x^\alpha \kappa(nx) = \infty$ (since $\kappa$ is slowly varying at infinity). The assumed continuity of $\kappa$ in $(0, \infty)$ does the rest.

Next we note that, for any integer $n \geq 1$,

$$(7.8) \qquad (n\varepsilon_n)^\alpha \kappa(n\varepsilon_n) = n^\alpha.$$

This implies that $\lim_{n \to \infty} n\varepsilon_n = +\infty$.

Now we estimate the integral in (7.7). It is easier to make a change of variables ($s := r/\varepsilon_n$) and deduce the following:

$$(7.9) \qquad \int_0^\infty e^{-r^\alpha \kappa(nr)} r^{\beta-1}\, dr = \varepsilon_n^\beta \int_0^\infty e^{-\varepsilon_n^\alpha \kappa(n\varepsilon_n s) s^\alpha} s^{\beta-1}\, ds := \varepsilon_n^\beta T_n.$$

Our goal is to show that $T_n \asymp 1$ $(n \to \infty)$. Note that if $\kappa(x) \asymp 1$ $(x \to \infty)$, then $\varepsilon_n \asymp 1$ and so $T_n \asymp 1$ $(n \to \infty)$.

In the general case, it is not a surprise that this is done by analyzing the integral over different regions; this is what we do next. We will need to make use of the representation theorem and the uniform convergence theorem for slowly varying functions; see Bingham, Goldie and Teugels (1987).

Thanks to (7.8), we have

$$(7.10) \qquad \begin{aligned} \int_1^\infty & e^{-\varepsilon_n^\alpha \kappa(n\varepsilon_n s) s^\alpha} s^{\beta-1}\, ds \\ &= \int_1^\infty \exp\!\left(-\varepsilon_n^\alpha \kappa(n\varepsilon_n) s^\alpha \frac{\kappa(n\varepsilon_n s)}{\kappa(n\varepsilon_n)}\right) s^{\beta-1}\, ds \\ &= \int_1^\infty \exp\!\left(-s^\alpha \frac{\kappa(n\varepsilon_n s)}{\kappa(n\varepsilon_n)}\right) s^{\beta-1}\, ds \\ &\leq \int_1^\infty \exp(-s^{\alpha-\delta}) s^{\beta-1}\, ds < \infty, \end{aligned}$$

where $0 < \delta < \alpha$ is a constant, and we have used the representation theorem for $\kappa$ in order to derive that

$$\frac{\kappa(n\varepsilon_n s)}{\kappa(n\varepsilon_n)} \geq s^{-\delta} \qquad \text{for all } n \text{ large enough.} \tag{7.11}$$

On the other hand, since $\kappa$ is nonnegative, we have

$$\int_0^1 e^{-\varepsilon_n^\alpha \kappa(n\varepsilon_n s)s^\alpha} s^{\beta-1}\,ds \leq \int_0^1 s^{\beta-1}\,ds = \frac{1}{\beta}. \tag{7.12}$$

Finally, it follows from (7.8) and the uniform convergence theorem for $\kappa$ that

$$\int_1^2 e^{-\varepsilon_n^\alpha \kappa(n\varepsilon_n s)s^\alpha} s^{\beta-1}\,ds = \int_1^2 \exp\left(-s^\alpha \frac{\kappa(n\varepsilon_n s)}{\kappa(n\varepsilon_n)}\right) s^{\beta-1}\,ds$$
$$\to \int_1^2 \exp(-s^\alpha) s^{\beta-1}\,ds \qquad \text{as } n \to \infty. \tag{7.13}$$

Combining (7.10), (7.12) and (7.13), we see that $T_n \asymp 1\ (n \to \infty)$, as asserted. □

LEMMA 7.5. *Under the conditions of Theorem 7.1, (7.5) holds.*

PROOF. Let $f(x) := x^\alpha \kappa(x)\ (x > 0)$. Because $f$ is regularly varying at infinity, it has an asymptotic inverse function $f^\leftarrow$ which is monotone increasing and regularly varying with index $1/\alpha$; see Bingham, Goldie and Teugels [(1987), page 28]. Furthermore, it follows from Proposition 1.5.15 of Bingham, Goldie and Teugels (1987) that $f^\leftarrow$ can be expressed as

$$f^\leftarrow(y) \sim y^{1/\alpha} \kappa^\#(y^{1/\alpha}) \qquad \text{as } y \to \infty, \tag{7.14}$$

where $\kappa^\#$ is the *de Bruijin conjugate* of $\kappa$.

Now we apply Lemma 7.4 with $n := |x|^{-1/\alpha}$, and recall (7.6), to deduce that $|x|^{\beta/\alpha} f_{d-\beta}(x) \asymp \varepsilon_n^\beta\ (|x| \to 0)$. For all $n \geq 1$, since $\varepsilon_n^\alpha \kappa(n\varepsilon_n) = 1$, we have $f(n\varepsilon_n) = n^\alpha$. Recall that $n\varepsilon_n \to \infty$ as $n \to \infty$, so our remarks on $f^\leftarrow$ prove that $\varepsilon_n \sim n^{-1} f^\leftarrow(n^\alpha) \sim \kappa^\#(n)\ (n \to \infty)$. Whence we have $|x|^{\beta/\alpha} f_{d-\beta}(x) \asymp [\kappa^\#(|x|^{-1/\alpha})]^\beta\ (|x| \to 0)$. This completes the proof of Lemma 7.5. □

REMARK 7.6. In order to prove Remark 7.3, we will use the following connection between $\kappa$ and its de Brujin conjugate $\kappa^\#$:

$$\kappa^\#(x) \sim [\kappa(x\kappa^\#(x))]^{-1} \qquad (x \to \infty); \tag{7.15}$$

see Bingham, Goldie and Teugels [(1987), Theorem 1.5.13]. Now we assume, in addition, that $\kappa(e^t) = t^\gamma \ell(t)$ for $t > 0$, where $\gamma$ is a constant and $\ell(\cdot)$ is slowly varying at infinity. Then we can write $\kappa(x) = (\ln x)^\gamma \ell(\ln x)$ for all $x > 1$. Consequently,

$$(7.16) \qquad \kappa(x\kappa^\#(x)) = \kappa(x)\left[1 + \frac{\ln \kappa^\#(x)}{\ln x}\right]^\gamma \frac{\ell(\ln x + \ln \kappa^\#(x))}{\ell(\ln x)}.$$

Since $\kappa^\#$ is slowly varying at infinity, we have $\ln \kappa^\#(x) = o(\ln x)$ as $x \to \infty$. This, and the representation theorem for $\ell(\cdot)$, together imply that $\ell(\ln x + \ln \kappa^\#(x)) \sim \ell(\ln x)$ as $x \to \infty$. Hence, it follows from (7.15) that

$$(7.17) \qquad \kappa^\#(x) \asymp \frac{1}{\kappa(x)} \qquad (x \to \infty).$$

Using again the assumption that $\kappa(e^t)$ is regularly varying at infinity, we deduce that Theorem 7.1 holds for the function $g_\kappa$ defined by (7.4).

7.2. *Dimension bounds: image.* For our next example, we consider the case where $X$ is an isotropic Lévy process in $\mathbf{R}^d$ and satisfies the following for two fixed constants $\delta, \eta \in (0, 2]$:

$$(7.18) \qquad \|\xi\|^{\delta+o(1)} \leq \Psi(\xi) \leq \|\xi\|^{\eta+o(1)} \qquad (\|\xi\| \to \infty).$$

A change of variables reveals that, for any $\beta \in (0, d)$,

$$(7.19) \qquad \frac{\beta}{\eta} \leq \liminf_{r \downarrow 1} \frac{\log f_{d-\beta}(r)}{\log(1/r)} \leq \limsup_{r \downarrow 0} \frac{\log f_{d-\beta}(r)}{\log(1/r)} \leq \frac{\beta}{\delta}.$$

Solve for the critical $\beta$ to see that $I(G) \geq \delta \dim G$ and $J(G) \leq \eta \dim G$. Thus, in this case,

$$(7.20) \qquad \delta \dim G \leq \dim X(G) \leq \eta \dim G \qquad \text{a.s.}$$

Note that the above includes the isotropic $\alpha$-stable processes, as well as Lévy processes with exponents that are regularly varying at infinity. Examples of the later processes can be found in Marcus (2001). More generally, a large class of Lévy processes satisfying (7.18) can be constructed by using the subordination method. Let $Y := \{Y(t)\}_{t \geq 0}$ be an isotropic $\alpha$-stable Lévy process in $\mathbf{R}^d$ and let $\tau := \{\tau(t)\}_{t \geq 0}$ be a subordinator with lower and upper indices $\sigma$ and $\beta$, respectively. Then the subordinated process $X := \{X(t)\}_{t \geq 0}$ defined by $X(t) := Y(\tau(t))$ is a Lévy process satisfying (7.18), with $\delta = \sigma \alpha$ and $\eta = \beta \alpha$. For other results along these lines see Blumenthal and Getoor (1961) and Millar (1971).

7.3. *Isotropic processes: preimage.* Suppose $X$ is isotropic, satisfies the absolute continuity condition of Corollary 3.2, and the regular variation condition (7.18) holds. Then, for any $\gamma \in (0, 1)$,

$$(7.21) \quad \|\xi\|^{\eta(\gamma-1)+o(1)} \leq \mathrm{Re}\left(\frac{1}{1+\Psi^{1-\gamma}(\xi)}\right) \leq \|\xi\|^{\delta(\gamma-1)+o(1)} \qquad (\|\xi\| \to \infty).$$

Now recall that, for any $\beta \in (0, d)$, the inverse Fourier transform of $\mathbf{R}^d \ni x \mapsto \|x\|^{-\beta}$ is a constant multiple of $\xi \mapsto \|\xi\|^{\beta-d}$. Thus, by Plancherel's theorem, for any $\nu \in \mathcal{P}(\mathbf{R}^d)$,

$$\mathcal{E}_\beta(\nu) = c_{d,\beta} \int_{\mathbf{R}^d} |\widehat{\nu}(\xi)|^2 \, \|\xi\|^{\beta-d} \, d\xi; \tag{7.22}$$

see Mattila [(1995), Lemma 12.12] and Kahane [(1985a), page 134]. Thus, thanks to the Frostman theorem (2.8), we have the following calculation in the present setting:

$$\frac{\eta + \dim R - d}{\eta} \leq \|\dim X^{-1}(R)\|_{L^\infty(\mathrm{P})} \leq \frac{\delta + \dim R - d}{\delta}. \tag{7.23}$$

When $\delta = \eta := \alpha$, (3.2) and (3.3) are ready consequences of this.

In fact, one can do more at little extra cost. Instead of isotropy, let us assume that $\Psi$ satisfies the *sector condition*: As $\|\xi\| \to \infty$, $\mathrm{Im}\Psi(\xi) = O(\mathrm{Re}\Psi(\xi))$. A few tedious, but routine, lines of calculations (see below) show that given any $\gamma \in (0, 1)$, $\Psi^{1-\gamma}$ also satisfies the sector condition, and so there exists a constant $c > 0$ such that, for all $\xi \in \mathbf{R}^d$,

$$\frac{c}{|1 + \Psi(\xi)|^{1-\gamma}} \leq \mathrm{Re}\left(\frac{1}{1 + \Psi^{1-\gamma}(\xi)}\right) \leq \frac{1}{|1 + \Psi(\xi)|^{1-\gamma}}. \tag{7.24}$$

If, in addition, there exist $\delta, \eta \in [0, 2]$ such that $\|\xi\|^{\delta+o(1)} \leq \mathrm{Re}\Psi(\xi) \leq \|\xi\|^{\eta+o(1)}$ ($\|\xi\| \to \infty$), then (7.23) holds. Another simple consequence of this example is that (3.4) continues to hold for all strictly $\alpha$-stable processes. We leave the details to the interested reader, and conclude this subsection by verifying the claim that whenever $\Psi$ satisfies the sector condition, then so does $\Psi^a$ for any $a \in \mathbf{R}$ with $|a| < 1$.

Write $\Psi(z) := |\Psi(z)|e^{i\theta(z)}$, where $\theta(z) \in [-\pi, \pi]$. By the sector condition on $\Psi$, there exists $c > 0$ such that, for all $\|\xi\|$ large enough, $|\mathrm{Im}\Psi(\xi)| \leq c\mathrm{Re}\Psi(\xi)$. But

$$|\sin(\theta(\xi))| = \frac{|\mathrm{Im}\Psi(\xi)|}{|\Psi(\xi)|} \leq \frac{c}{\sqrt{1+c^2}} := \sin(\eta) < 1, \tag{7.25}$$

where $\eta := \sin^{-1}(c/\sqrt{1+c^2})$. This means that for any fixed $a \in \mathbf{R}$ with $|a| < 1$, $\cos(a\theta(\xi)) \geq \cos(a\eta) > 0$ as soon as $\|\xi\|$ is large enough. Therefore, there exists $\varepsilon := \cos(a\eta) > 0$ such that for any $a \in \mathbf{R}$ with $|a| < 1$, and all $\|\xi\|$ large,

$$\mathrm{Re}\Psi^a(\xi) = |\Psi(\xi)|^a \cos(a\theta(\xi)) \geq \varepsilon|\Psi(\xi)|^a = \varepsilon|\Psi^a(\xi)| \geq |\mathrm{Im}(\Psi^a(\xi))|. \tag{7.26}$$

This proves that the sector condition holds for $\Psi^a$.

7.4. *Processes with stable components.* A (Lévy) *process $X$ with stable components* is a $d$-dimensional process with independent components $X_1, \ldots, X_p$ such that $X_j$ is an $\alpha_j$-stable Lévy process in $\mathbf{R}^{d_j}$, where $d = \sum_{j=1}^{p} d_j$. By relabelling the components, we can and will assume throughout that $2 \geq \alpha_1 \geq \alpha_2 \geq \cdots \geq \alpha_p > 0$.

Pruitt and Taylor (1969) have studied the range of $X$, and proved that, with probability one,

$$(7.27) \qquad \dim X(\mathbf{R}_+) = \begin{cases} \alpha_1, & \text{if } \alpha_1 \leq d_1, \\ 1 + \alpha_2(1 - \alpha_1^{-1}), & \text{if } \alpha_1 > d_1 = 1. \end{cases}$$

Becker-Kern, Meerschaert and Scheffler (2003) have recently extended (7.27) to a class of operator-stable Lévy processes in $\mathbf{R}^d$, which allow dependence among the components $X_1, \ldots, X_p$. Their argument involves making a number of technical probability estimates, and makes heavy use of the results of Pruitt (1969). As a result, they impose some restrictions on the transition densities of $X$.

In the following, we give a different analytic proof of the result (7.27). Since we do not need probability estimates, our argument works for more general Lévy processes than those of Pruitt and Taylor (1969). In particular, we expect that our method will work for the cases that have remained unsolved by Becker-Kern, Meerschaert and Scheffler (2003).

PROPOSITION 7.7. *Let $X$ be a Lévy process in $\mathbf{R}^d$, with $d \geq 2$, whose Lévy exponent $\Psi$ satisfies the following:*

$$(7.28) \qquad \mathrm{Re}\left(\frac{1}{1+\Psi(\xi)}\right) \asymp \frac{1}{\sum_{j=1}^{p}|\xi_j|^{\alpha_j}} \qquad as \ \|\xi\| \to \infty.$$

*Then almost surely,*

$$(7.29) \qquad \dim X(\mathbf{R}_+) = \begin{cases} \alpha_1, & \text{if } \alpha_1 \leq d_1, \\ 1 + \alpha_2(1 - \alpha_1^{-1}), & \text{if } \alpha_1 > d_1. \end{cases}$$

REMARK 7.8. Condition (7.28) is satisfied by a large class of Lévy processes, including the Lévy processes with stable components considered by Pruitt and Taylor (1969), as well as more general operator-stable Lévy processes. Moreover, one can replace the power functions $|\xi_j|^{\alpha_j}$ by regularly varying functions and the conclusion still holds. In particular, (7.29) still holds if $X$ is a Lévy process in $\mathbf{R}^d$ whose components involve independent asymmetric Cauchy processes.

PROOF OF PROPOSITION 7.7. For any $\gamma > 0$, it follows from (7.28) that the integral in (4.18) is comparable to

$$(7.30) \qquad I_\gamma := \int_{\{\xi \in \mathbf{R}^d : \|\xi\| \geq 1\}} \frac{1}{1 + \sum_{j=1}^{p}|\xi_j|^{\alpha_j}} \frac{d\xi}{\|\xi\|^{d-\gamma}}.$$

Notice that $I_\gamma = \infty$ for all $\gamma \geq \alpha_1$. Hence, we always have $\dim X(\mathbf{R}_+) \leq \alpha_1$ almost surely (Corollary 4.7).

Now we derive the corresponding lower bound in the case that $\alpha_1 \leq d_1$. It is sufficient to work with the two-dimensional Lévy process $X = (X_1, X_2)$. Hence, without loss of generality, we will assume that $d = 2$.

Clearly, if $d_1 = d = 2$, then it follows from (7.30) that $I_\gamma < \infty$ for all $0 < \gamma < \alpha_1$. Thus, Corollary 4.7 implies $\dim X(\mathbf{R}_+) \geq \alpha_1$ almost surely, as desired. So we only need to consider the case when $d_1 = 1$ and $\alpha_1 \leq 1$. Write

$$
\begin{aligned}
I_\gamma &\asymp \left[\int_0^1 \int_1^\infty + \int_1^\infty \int_0^1\right] \frac{1}{1+\xi_1^{\alpha_1}+\xi_2^{\alpha_2}} \cdot \frac{d\xi_1\,d\xi_2}{\xi_1^{2-\gamma}+\xi_2^{2-\gamma}} \\
&\quad + \int_1^\infty d\xi_1 \int_1^\infty \frac{1}{1+\xi_1^{\alpha_1}+\xi_2^{\alpha_2}} \cdot \frac{d\xi_2}{\xi_1^{2-\gamma}+\xi_2^{2-\gamma}} \\
&:= I_\gamma^{(1)} + I_\gamma^{(2)}.
\end{aligned}
\tag{7.31}
$$

For any $0 < \gamma < \alpha_1 \leq 1$, $I_\gamma^{(1)}$ is finite, and

$$
\begin{aligned}
I_\gamma^{(2)} &\leq \int_1^\infty \frac{d\xi_1}{1+\xi_1^{\alpha_1}} \cdot \int_1^\infty \frac{d\xi_2}{\xi_1^{2-\gamma}+\xi_2^{2-\gamma}} \\
&\leq \int_1^\infty \frac{1}{1+\xi_1^{\alpha_1}} \cdot \frac{d\xi_1}{\xi_1^{1-\gamma}} \cdot \int_0^\infty \frac{d\xi_2}{1+\xi_2^{2-\gamma}} < \infty.
\end{aligned}
\tag{7.32}
$$

Consequently, $I_\gamma < \infty$ for all $\gamma < \alpha_1$. It follows from Corollary 4.7 that, when $\alpha_1 \leq d_1$, $\dim X(\mathbf{R}_+) \geq \alpha_1$ almost surely. This proves the first part of (7.29).

Next we prove the second part of (7.29). Since $\alpha_1 > d_1 = 1$, we have $\alpha_2 \leq 1 + \alpha_2(1 - \alpha_1^{-1}) \leq \alpha_1$. For any $\gamma > 1 + \alpha_2(1 - \alpha_1^{-1})$, in order to prove that $I_\gamma = \infty$, we will make use of the following inequality: If $d > 1 + \gamma$ and $\alpha > 0$, then for all constants $a, b \geq 2$ that satisfy $b^{1/\alpha} a^{-1} \geq K_1^{-1}$,

$$
\begin{aligned}
\int_1^\infty \frac{1}{b+x^\alpha} &\cdot \frac{dx}{(a^2+x^2)^{(d-\gamma)/2}} \\
&= a^{-(d-1-\gamma)} \int_{a^{-1}}^\infty \frac{1}{b+a^\alpha x^\alpha} \cdot \frac{dx}{(1+x^2)^{(d-\gamma)/2}} \\
&\geq a^{-(d-1-\gamma)} \int_{a^{-1}}^{K_1 b^{1/\alpha} a^{-1}} \frac{1}{b+a^\alpha x^\alpha} \cdot \frac{dx}{(1+x^2)^{(d-\gamma)/2}} \\
&\geq K_2 b^{-1} a^{-(d-1-\gamma)},
\end{aligned}
\tag{7.33}
$$

where $K_1$ and $K_2$ are positive and finite constants.

We rewrite the integral in (7.30) in all $d$ coordinates and relabel $\alpha_1, \ldots, \alpha_p$ for each coordinate in an obvious way (now denoted as $\alpha_1, \ldots, \alpha_d$) to derive

$$
I_\gamma \geq \int_1^\infty d\xi_1 \cdots \int_1^\infty \frac{1}{1+\xi_1^{\alpha_1}+\cdots+\xi_d^{\alpha_d}} \cdot \frac{d\xi_d}{\|\xi\|^{d-\gamma}}.
\tag{7.34}
$$

If $d > 2$, then we iteratively integrate the integral in (7.34) $[d\xi_d \times d\xi_{d-1} \times \cdots \times d\xi_3]$, and use (7.33) $d-2$ times. (Note that, for the obvious choices of $a$ and $b$, the condition $b^{1/\alpha}a^{-1} \geq K_1^{-1}$ holds for some constant $K_1 > 0$ because of the assumption $\alpha_1 \geq \alpha_2 \geq \cdots \geq \alpha_d$.) As a result, we deduce that there is a constant $K_3 > 0$ such that

$$(7.35) \qquad I_\gamma \geq K_3 \int_1^\infty d\xi_1 \int_1^\infty \frac{1}{\xi_1^{\alpha_1} + \xi_2^{\alpha_2}} \cdot \frac{d\xi_2}{\xi_1^{2-\gamma} + \xi_2^{2-\gamma}} := K_3 I_\gamma^{(3)}.$$

Clearly, this inequality also holds for $d = 2$. A change of variables then yields

$$(7.36) \quad \begin{aligned} I_\gamma^{(3)} &= \int_1^\infty \frac{d\xi_1}{\xi_1^{1+\alpha_2-\gamma}} \int_{\xi_1^{-1}}^\infty \frac{1}{\xi_1^{\alpha_1-\alpha_2} + x^{\alpha_2}} \cdot \frac{dx}{1 + x^{2-\gamma}} \\ &\geq \frac{1}{2} \int_1^\infty \frac{d\xi_1}{\xi_1^{1+\alpha_2-\gamma}} \int_1^\infty \frac{1}{\xi_1^{\alpha_1-\alpha_2} + x^{\alpha_2}} \cdot \frac{dx}{x^{2-\gamma}} \\ &\geq \frac{1}{2} \int_1^\infty \frac{d\xi_1}{\xi_1^{1+\alpha_2-\gamma}} \cdot \xi_1^{-(\alpha_1-\alpha_2)(1+(1-\gamma)/\alpha_2)} \int_1^\infty \frac{1}{1+y^{\alpha_2}} \cdot \frac{dy}{y^{2-\gamma}} \\ &\geq K_4 \int_1^\infty \frac{d\xi_1}{\xi_1^{\alpha_1+(1-\gamma)\alpha_1/\alpha_2}}. \end{aligned}$$

Recall that $\gamma > 1 + \alpha_2(1 - \alpha_1^{-1})$. Equivalently, we have $\alpha_1 + (1-\gamma)\alpha_1/\alpha_2 \leq 1$. Combining (7.34)–(7.36) together yields $I_\gamma = \infty$; this proves that $\dim X(\mathbf{R}_+) \leq 1 + \alpha_2(1 + \alpha_1^{-1})$, a.s. (Corollary 4.7).

Finally, we prove the lower bound for $\dim X(\mathbf{R}_+)$ in the case that $\alpha_1 > d_1 = 1$. Again, it suffices to assume that $d = 2$; otherwise, consider the projection of $X$ into $\mathbf{R}^2$. For any $1 < \gamma < 1 + \alpha_2(1 - \alpha_1^{-1})$, we have $2 - \gamma + \alpha_2 > 1$; hence, (7.31) implies that there exist positive and finite constants $K_5$ and $K_6$ such that

$$(7.37) \qquad I_\gamma \leq K_5 + K_6 I_\gamma^{(3)}.$$

As we did for (7.36), we can prove that

$$(7.38) \quad \begin{aligned} I_\gamma^{(3)} &= \int_1^\infty \frac{d\xi_1}{\xi_1^{1+\alpha_2-\gamma}} \left[ \int_{\xi_1^{-1}}^1 + \int_1^\infty \right] \frac{1}{\xi_1^{\alpha_1-\alpha_2} + x^{\alpha_2}} \cdot \frac{dx}{1+x^{2-\gamma}} \\ &\leq \int_1^\infty \frac{d\xi_1}{\xi_1^{1+\alpha_1-\gamma}} + \int_1^\infty \frac{d\xi_1}{\xi_1^{1+\alpha_2-\gamma}} \int_1^\infty \frac{1}{\xi_1^{\alpha_1-\alpha_2} + x^{\alpha_2}} \cdot \frac{dx}{1+x^{2-\gamma}} \\ &\leq \int_1^\infty \frac{d\xi_1}{\xi_1^{1+\alpha_1-\gamma}} + \int_1^\infty \frac{d\xi_1}{\xi_1^{\alpha_1+(1-\gamma)\alpha_1/\alpha_2}} \int_0^\infty \frac{1}{1+x^{\alpha_2}} \cdot \frac{dx}{x^{2-\gamma}}. \end{aligned}$$

Observe that all three integrals in (7.38) are finite because $1 < \gamma < 1 + \alpha_2(1 - \alpha_1^{-1}) < \alpha_1$ and $2 - \gamma + \alpha_2 > 1$. It follows from (7.37) that $I_\gamma < \infty$ for all $\gamma < 1 + \alpha_2(1 - \alpha_1^{-1})$. Hence, Corollary 4.7 implies that $\dim X(\mathbf{R}_+) \geq 1 + \alpha_2(1 - \alpha_1^{-1})$, a.s. This finishes the proof of Proposition 7.7. □

**Acknowledgment.** We are grateful to an anonymous referee for his/her careful reading, as well as for pointing out some misprints in the first draft of this paper.

Department of Mathematics
University of Utah
155 South 1400 East JWB 233
Salt Lake City, Utah 84112–0090
USA
e-mail: davar@math.utah.edu
url: www.math.utah.edu/~davar

Department of Statistics
and Probability
Michigan State University
A-413 Wells Hall
East Lansing, Michigan 48824
USA
e-mail: xiao@stt.msu.edu
url: www.stt.msu.edu/~xiaoyimi